\documentclass[11pt,twoside]{article}

\usepackage{achemso}
\setkeys{acs}{articletitle=true}
\usepackage{epsfig,amsfonts,xcolor}
\usepackage{amsmath}
\usepackage{lineno}
\usepackage{pdflscape}
\usepackage{subfigure}

\usepackage{amssymb, palatino,geometry,url}
\usepackage{algorithmic}
\usepackage[colorlinks=false,linkcolor=black,citecolor=black,urlcolor=black]{hyperref}
\usepackage{soul}
\allowdisplaybreaks

\usepackage{fancyhdr}
\begin{document}

\title{\Large Graph-Based Optimization for Technology Pathway Analysis:\\  A Case Study in Decarbonization of University Campuses}

\author{Blake Lopez${}^{\dag}$, Jiaze Ma${}^{\dag}$, Victor M. Zavala${}^{\dag}$\thanks{Corresponding Author: victor.zavala@wisc.edu.}\\
 \\
  {\small ${}^{\dag}$Department of Chemical and Biological Engineering}\\
 {\small \;University of Wisconsin - Madison, 1415 Engineering Dr, Madison, WI 53706, USA}}

 \date{}
\maketitle
\begin{abstract}
Industrial sectors such as urban centers, chemical companies, manufacturing facilities, and microgrids are actively exploring strategies to help reduce their carbon footprint. For instance, university campuses are complex urban districts (involving collections of buildings and utility systems) that are seeking to reduce carbon footprints that originate from diverse activities (e.g., transportation operations and production of heating, cooling, and power utilities). This work presents an optimization framework to identify technology pathways that enable decarbonization of complex industrial sectors. The framework uses a graph abstraction that compactly captures interdependencies between diverse products and technologies as well as diverse externalities (e.g., market, policy, and carbon prices). Duality analysis reveals that the formulation can be interpreted as an economy, market, or value chain that uses technologies to generate economic value (wealth) by transforming basic products into higher value products. This interpretation also reveals that the formulation identifies pathways that maximize the profit of stakeholders, helps reveal the inherent value (prices) of intermediate products, and helps analyze the impact of externalities and technology specifications on product values. Our developments are illustrated via a case study involving a prototypical university campus that seeks to identify pathways that reduce its carbon footprint (e.g., via electrification and deployment of hydrogen technologies). We use the framework to determine carbon tax values, technology specifications, and investment budgets that activate different technology pathways and that achieve different levels of decarbonization.
\end{abstract}


\newpage
%%%%%%%%%%%%%%%%%%%%%%%%%%%%%%%%%%%%%%%%%%

\section{Introduction}

The global mean surface temperature reached 1℃ above pre-industrial times in 2017 and is projected to reach 1.5℃ as early as 2030 \cite{masson-delmotte_ipcc_nodate}. Significant efforts need to be made to decrease greenhouse gas emissions to limit warming. The most notable effort to accomplish this is the Paris agreement, which seeks to limit warming to at most 2℃ by providing technological and financial support as well as capacity building between the 192 countries in the agreement \cite{conference_of_the_parties_to_the_united_nations_framework_convention_on_climate_change_21st_sess__2015__paris_report_2016}. To satisfy the goals of global agreements, industrial sectors such as chemical companies, manufacturing facilities, farms, microgrids, and urban centers need to shift their operations so as to minimize CO\textsubscript{2} emissions. For instance, ammonia production facilities are exploring the implementation of electrolysis technologies that harness renewable power to replace hydrogen sourced from natural gas, and utility companies are relying on falling renewable costs and energy storage technologies to mitigate carbon emissions \cite{luderer_impact_2022,tullo_industry_2021,tullo_cf_2020}. 
\\

An important and representative set of industrial systems that is actively seeking to decarbonize operations are {\em university campuses}. To give some perspective into the magnitude of these systems, we note that there are over 16.5 million students in the US reside on university campuses \cite{noauthor_nces_nodate}. Campuses involve large collections of buildings that consume heating, cooling, power, and transportation services. Many university campuses own and operate utility systems that provide steam, hot water, and chilled water to buildings and research facilities. Currently, most of these systems are powered using combined heat and power (CHP) systems that run on natural gas \cite{sinha_greenhouse_2010}. Electrical power is also used for diverse essential services (e.g., ventilation, lighting, circulation, lectures) and for resource-intensive equipment such as computational facilities (e.g., computing clusters and data centers). The generation and distribution of all these utilities generate a massive carbon footprint. Moreover, universities also typically operate transportation services, which further contribute to the carbon footprint. In essence, university campuses are representative systems in that they combine elements of urban centers and manufacturing facilities. It has been estimated that university campuses in the US emit 7.7 MTCO\textsubscript{2}e per student and are responsible for 2\% of the 5,222 million MTCO\textsubscript{2}e of the US greenhouse gas emissions \cite{sinha_greenhouse_2010,us_epa_inventory_2017}. Notably, these emissions are higher than those of ammonia production.
\\

Technologies such as renewable power generation, heat pumps, electric and hydrogen-based transportation vehicles, and energy storage systems are actively being studied as means to help decrease carbon footprints \cite{coppez_importance_2010,lah_decarbonizing_2017,vargas_wind_2019,roberts_current_2016,khan_solar_2016}. The number of new technologies that can be used is steadily increasing and complicates decision-making processes, because it is important to understand how new technologies will interact with existing infrastructure and operations. In the context of university campuses, cost-minimization models have been proposed to determine the energy mix, capacities, and storage requirements for a university campus (model is cast as a mixed integer nonlinear program) \cite{tian_sustainable_2022}. This modeling approach requires detailed data for investment costs, operating costs, operating lifetimes, and scaling factors for the technologies proposed and time-dependent demand data to determine a minimized cost scenario with no emissions. Other modeling approaches have used multiobjective optimization to determine technology operations that can help reach emissions reduction targets, determine lowest cost for maximum renewable use, and determine cost of reducing emissions \cite{olsen_planning_2019,mirzaei_alavijeh_cost-effectiveness_2020,ho_multi-objective_2014}. Existing optimization models are typically intended to provide long-term investment plans to achieve decarbonization. Besides models seeking to decarbonize, technology selection and the general exploration of pathways can be seen applied frequently in the conversion of a raw feedstock or waste product to more valuable products \cite{corsano_optimal_2011,bao_shortcut_2011, morsy_superstructure_2023}. These approaches often take the approach of a cost minimization but can also seek to maximize the yield or a measure of the economic potential or profit. The models in case studies tend to be case specific and require detailed technology specification. Detailed models are ought to be used to make important investment decisions but these often hindered by the fact that there might not be sufficient data to evaluate emerging technologies; in addition, the use of detailed models might require significant implementation and computational effort. As such, there is also a need for simple models that can help screen options prior to investing resources in model detailed studies. 
\\

In this work, we provide an optimization framework for analyzing technology pathways that can provide desired product demand targets while minimizing supply and technology costs. The framework uses a compact, graph/network representation that aims to facilitate the analysis of complex interdependencies that arise between products and technologies. In the proposed representation, we consider the supply of basic products to the system (e.g., natural gas and electricity) and technologies that transform these basic products to higher-value intermediate and final products (e.g., hydrogen and cooling water) that are delivered to consumers. The formulation operates at a high level of abstraction, which enables capturing diverse types of products (e.g., material, energy, labor, services) and byproducts (e.g., carbon emissions and waste) in a unified manner. The proposed formulation can be viewed as a superstructure network that interconnects different stakeholders in a system and helps determine how externalities (e.g., policy and technology specifications) impact the system. Moreover, duality analysis reveals that the formulation has a natural economic/market interpretation; this reveals that the formulation aims to identify technology pathways that maximize total economic value (by using technologies for transforming products into higher-value products) and explains how externalities can help activate/deactive pathways. The market interpretation also allows us to discover inherent values (prices) for key intermediate products and helps evaluate the impact of externalities and technology specifications on such values. Moreover, the framework provides insights into how expanding the domain of the system (e.g., by considering alternative products) can activate technology pathways. The modeling framework uses minimal data on products and technologies, which facilitates analysis and enables fast screening of diverse technologies and externalities. We provide a case study that analyzes decarbonization pathways for a prototypical university campus to illustrate the developments. 
\\

%%%%%%%%%%%%%%%%%%%%%%%%%%%%%%%%%%%%%%%%%%

\section{Graph-Based Optimization Framework}

We consider a system comprised of a set of basic, intermediate, and final products as well as a set of technologies that transform basic products into intermediate and final products. The goal is to determine technology pathways that can maximize the value of served product demands, while minimizing supply and technology costs. 
\\

A couple of optimization models are introduced to determine optimal technology pathways that conduct desired goals. A {\em management model} will be used to understand how suppliers, technologies, and demands interact and to determine optimal allocations for all these stakeholders. The model uses a graph abstraction that captures the topology/connectivity that arises from interactions between products and technologies. The model can be interpreted as a value chain that aims to identify pathways that maximize profit for all stakeholders involved; this is done by using technologies to create wealth (transforming lower-value products into higher-value products). This interpretation will also allows us to determine the inherent value (prices) of products, which is particularly useful in attributing value for intermediate products and to understand how externalities (e.g., carbon prices or disposal costs for waste products) and technology costs/efficiencies propagate through the technology pathways and influence product prices. An {\em investment model} will be used to prioritize the selection of pathways under constrained investment budgets and to trade-off profits with investment costs.
\\

Our models aim to use minimum specification data for technologies (e.g., efficiencies, capacities, operating costs, investment costs), so as to provide $``$high-level$"$ picture of the technology landscape and on potential factors that influence their economic viability. This is also, in part, motivated by the fact that there is often limited data available for existing technologies. As expected, any simplification to the representation of technologies will come with a downside of inaccuracy. Our high-level abstraction aims to provide an intuitive approach that helps navigate complex interdependencies between products,  technologies, and externalities. 
\\

The proposed model can be interpreted as a value chain in which there is no transportation and spatial/geographical context; specifically, we note that our model is a simplification of the multi-product supply chain model presented in \cite{sampat_coordinated_2019,tominac_economic_2021}. This observation is important, as the model proposed here can provide valuable preliminary insights that can inform the development of more sophisticated models; for instance, it can help determine technology pathways that will or will not be selected in supply chain designs. 
\\

\subsection{Management Model}

We define a {\em system} comprised of a set of products $\mathcal{P}$, set of suppliers $\mathcal{S}$ that offer products, technologies $\mathcal{T}$ that transform products, and consumers $\mathcal{D}$ that request products. These elements are interconnected via a {\em graph} that captures product-technology dependencies. We can think of this system as an economy, market, or value chain that aims to generate economic value (wealth) by transforming basic products into higher value products. In other words, we think about this system as an economy in which elements (suppliers, consumers, technologies) are stakeholders that provide/request services and that aim to generate wealth. Factors that affect the ability of this economy from generating wealth include costs of basic products, technology factors (e.g., costs, capacities, and efficiencies), and externalities (e.g., policy, external markets, and taxes). 
\\

Each supplier $i\in\mathcal{S}$ has an offered value $\alpha_i\in \mathbb{R}$, associated allocation $s_i\in \mathbb{R}$ (flow of product provided to the system), an available capacity $\bar{s}_i\in \mathbb{R}$ (maximum product amount that it can provide), and a supplied product $p(i)\in \mathcal{P}$. Each consumer $j\in \mathcal{D}$ has an associated offered value $\alpha_j\in \mathbb{R}$, allocation $d_j\in\mathbb{R}$ (flow of product extracted from the system), an available capacity $\bar{d}_j\in \mathbb{R}$ (maximum product amount that it can take), and a requested product $p(j)\in \mathcal{P}$. We denote the set of suppliers or consumers that provide/request a given product $p\in\mathcal{P}$ as $\mathcal{S}_p\subseteq \mathcal{S}$ and $\mathcal{D}_p\subseteq\mathcal{D}$, respectively.
\\

Each technology $k\in \mathcal{T}$ has an offered value $\alpha_k\in \mathbb{R}$ (service cost), has an allocation flow $t_k\in \mathbb{R}$ (flow processed by the technology); this flow corresponds to a reference (input) product $p(k)\in \mathcal{P}$ and an available processing capacity $\bar{t}_k\in \mathbb{R}$ (maximum amount of reference product that it can process). Technologies are key assets, as they generate wealth by conducting transformation of products; this is captured using the transformation factors $\gamma_{k,p}\in \mathbb{R},\; k\in \mathcal{T},p\in\mathcal{P}$. A positive transformation factor ($\gamma_{k,p}>0$) indicates that product $p$ is produced by technology $k$, a negative factor ($\gamma_{k,p}<0$) indicated that the product is consumed by the technology, and a zero factor $\gamma_{k,p}=0$ indicates that the product is neither produced nor consumed (does not participate in the technology). The set of technologies that generate or consume a specific product $p\in \mathcal{P}$ are denoted as $\mathcal{T}_{p}\subseteq \mathcal{T}$. Transformation factors play a key role in dictating the behavior of the system, as they provide measures of efficiency for technologies; moreover, these factors capture product-technology interdependence (connectivity) and capture the generation of desirable (e.g., value-added) products and undesirable (e.g., waste) byproducts. In other words, transformation of products that occur in technologies define the topology of the system graph; this topology captures diverse pathways that might exist to obtain specific products and encodes key information that explains how diverse externalities propagate through the system. 
\\

\begin{figure}[!htp]
\center{\includegraphics[width=0.7\textwidth]{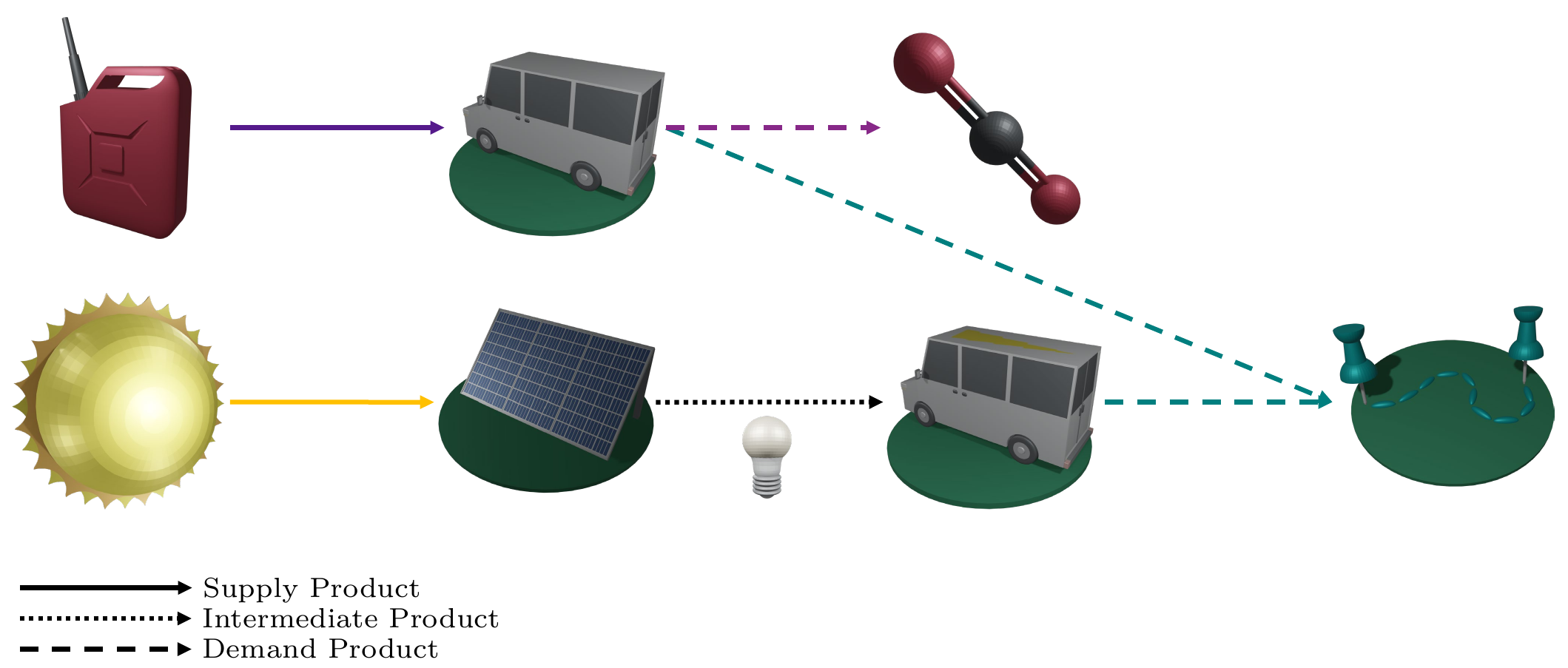}}
\caption{Illustration of a simple graph in which supplies of basic products (gasoline and solar radiation) are used by technologies (gas-powered vehicles,  electricity vehicles, and solar panels) to satisfy demand of higher value products (travel distance) and that generate undesired products (CO${}_2$ emissions). In this case, electrical power generated by the photovoltaic panels is an intermediate product that is used by the electric vehicle to generate a value-added and final product (distance). If the cost of CO${}_2$ emissions is not taken into account (and/or gasoline is inexpensive), the gas vehicle pathway will be preferred. On the other hand, when emissions costs are taken into account and/or renewable power is inexpensive, the electric vehicle pathway will be preferred.}
\label{fig:cars}
\end{figure}

\begin{figure}[!htp]
\center{\includegraphics[width=1.1\textwidth]{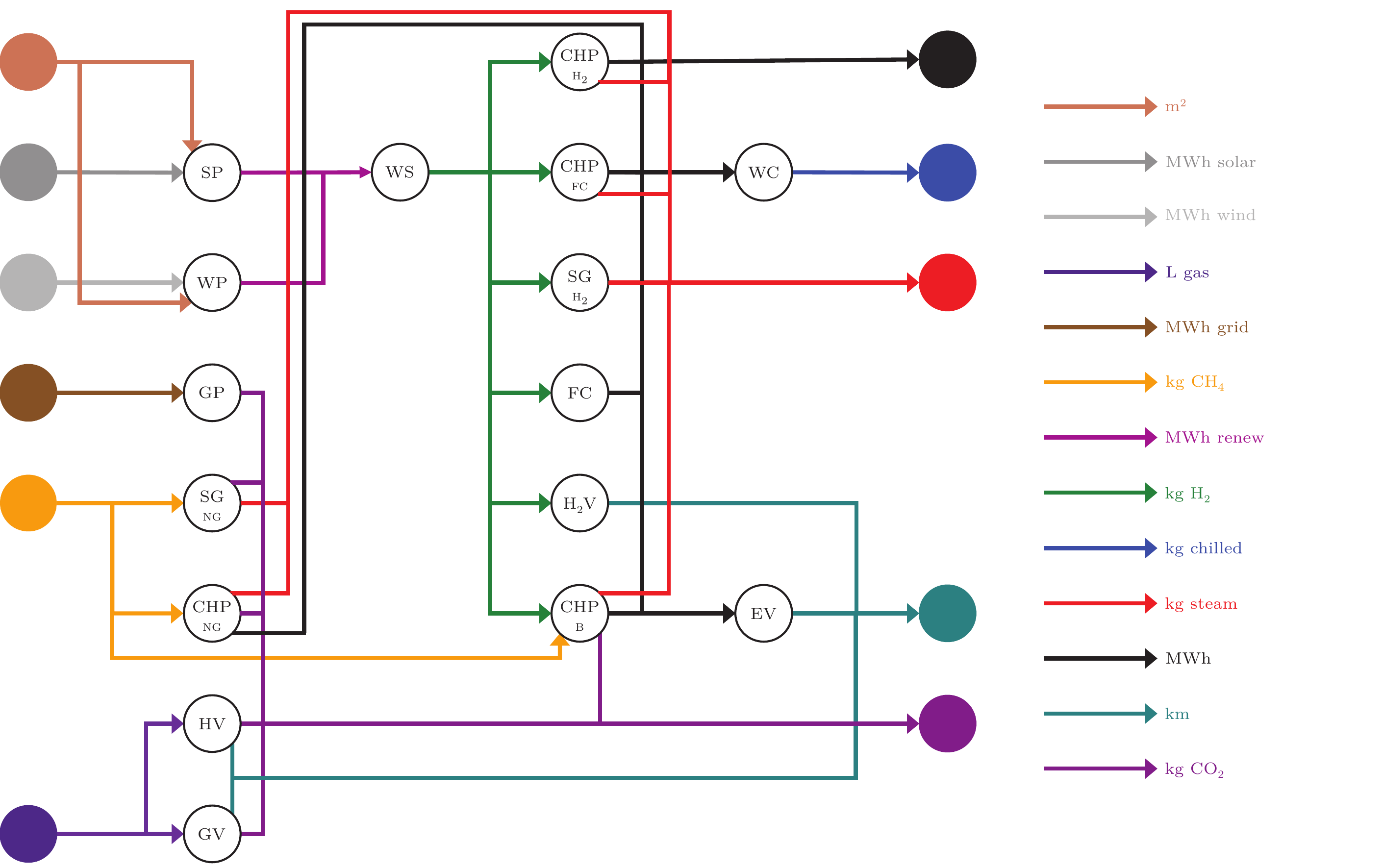}}
\caption{Graph showing diverse suppliers, technologies, and consumers that participate in a university campus. Complex interdependencies/connectivity arises between products and technologies and it is thus nonobvious what are the best pathways. See Table \ref{tab:Acronyms} for technology descriptions.}
\label{fig:SS}
\end{figure}

The {\em data} for suppliers, consumers, and technologies is summarized into a graph (a superstructure) that captures all product-technology interdependencies. The goal is to determine pathways in this superstructure that generate maximum wealth, while filtering out those that do not generate wealth. Figure \ref{fig:cars} provides a simple illustration on how elements are connected in the graph and how competing technology pathways emerge. Moreover, this aims to illustrate how externalities (emissions costs/taxes) can favor some pathways over others. A more complex graph showing diverse pathways to decarbonize a university campus is presented in Figure \ref{fig:SS}. This illustrates how complexity quickly arises due to the presence of diverse pathways and product-technologies connectivity; as such, it is necessary to use systematic techniques to identify optimal pathways. 
\\

We now propose an optimization formulation that helps identify optimal technology pathways from the graph. The formulation is given by: 
\begin{subequations} 
\label{MPrimal}
\begin{align} 
\label{MObj}
\max_{(s,d,t)} &\quad \sum_{j \in \mathcal{D}}{\alpha_{j}}d_{j} - \sum_{i \in \mathcal{S}}{\alpha_{i}}s_{i} - \sum_{k\in\mathcal{T}}{\alpha_{k}t_{k}}\\
\label{MBalance}
\text{s.t} &\quad \sum_{i \in \mathcal{S}_p}{s_i} + \sum_{k \in \mathcal{T}_{p}}{\gamma_{k,p}t_{k}} - \sum_{j \in \mathcal{D}_{p}}{d_{j}}=0, \quad p \in \mathcal{P}, \quad (\pi_{p})\\
\label{MS}
&\quad 0 \le s_{i} \le \bar{s}_{i}, \quad i \in \mathcal{S}\\
\label{MD}
&\quad 0 \le d_{j} \le \bar{d}_{j}, \quad j \in \mathcal{D}\\
\label{MT}
&\quad 0 \le t_{k} \le \bar{t}_{k}, \quad k \in \mathcal{T}
\end{align}
\end{subequations}

This $``$management$"$ model aims to identify product allocations $(s,d,t)$ for all system elements that maximize the total surplus \eqref{MBalance}. The total surplus balances the value of demand served (to be maximized) and the costs of all supplies and technologies (to be minimized). We will see (via duality analysis) that the total surplus is equivalent to the total system profit (revenue minus cost) and captures the profit of all system elements. 
\\

The constraints \eqref{MBalance} encode the graph that captures the connectivity between products and technologies. In this graph, products are interpreted as {\em nodes} and we enforce conservation of product at such nodes. Specifically, for each product $p\in \mathcal{P}$, we must have that all input by suppliers, generation/consumption by technologies, and extraction by consumers must be balanced. 
\\

The dual variables $\pi_p, p \in \mathcal{P}$ of the balance constraints \eqref{MBalance} can be interpreted as inherent values (prices) for the products. Specifically, the dual variables capture how valuable a given product is in the economy. For instance, a high-value product can be critical because it enables  diverse technology pathways and generate wealth via sales of final products; conversely, a low-value product can be of low priority (or even irrelevant) for the generation of wealth in the system. In fact, we will see later that prices can have negative values, indicating that certain products (e.g., waste) can be detrimental to the creation of wealth in the economy. However,  products with negative values might be necessary to achieve alternate goals in the system, such as mitigation of social and environmental impacts by the economy. The value of a product ($\pi_p$) can be different from the supply ($\alpha_i$) and demand ($\alpha_j$) bids because the bids are set by the suppliers and consumers while the value is determined by the dual formulation. However, the product values do depend on bids which will be shown later using the dual formulation of the problem. Determining prices for intermediate products is particularly important, as such products typically do not have external markets (they are means to an end).
\\

We can gain important insights into the economic properties of the optimization problem by analyzing the (partial) Lagrangian dual function:
\begin{equation} 
\label{ML}
\begin{aligned}
\mathcal{L}(s,d,t,\pi) = \sum_{j \in \mathcal{D}}{\alpha_{j}}d_{j} - \sum_{i \in \mathcal{S}}{\alpha_{i}}s_{i} - \sum_{k\in\mathcal{T}}{\alpha_{k}t_{k}} + \sum_{p \in \mathcal{P}}\pi_{p}\left ({\sum_{i \in \mathcal{S}_p}{s_i} + \sum_{k \in \mathcal{T}_{p}}{\gamma_{k,p}t_{k}} - \sum_{j \in \mathcal{D}_{p}}{d_{j}}}\right).
\end{aligned}
\end{equation}
If we now consider the following identities:
\begin{subequations} 
\label{ID}
\begin{align} 
\label{IDS}
\sum_{p \in \mathcal{P}}\pi_{p}\sum_{i \in \mathcal{S}_{p}}{s_i}= \sum_{i \in \mathcal{S}}{\pi_{i}s_{i}}\\
\label{IDT}
\sum_{p \in \mathcal{P}}\pi_{p}\sum_{k \in \mathcal{T}_{p}}{\gamma_{k,p}t_{k}} = \sum_{k \in \mathcal{T}}{\pi_{k}t_{k}}\\
\label{IDD}
\sum_{p \in \mathcal{P}}\pi_{p}\sum_{j \in \mathcal{D}_{p}}{d_j} = \sum_{j \in \mathcal{D}}{\pi_{j}d_{j}}, 
\end{align}
\end{subequations}
we can see that the Lagrange function can be rewritten as:
\begin{equation} 
\label{MLS}
\begin{aligned}
\mathcal{L}(s,d,t,\pi) = \sum_{j \in \mathcal{D}}{(\alpha_{j}^{\mathcal{D}} - \pi_j)d_{j}} + \sum_{i \in \mathcal{S}}{(\pi_i - \alpha_{i}^{\mathcal{S}})s_{i}} + \sum_{k \in \mathcal{T}}{(\pi_k - \alpha^{\mathcal{T}}_{k})t_{k}}.
\end{aligned}
\end{equation}
Here, we use the short-hand notation $\pi_i:=\pi_{p(i)}$ to denote the price of the product that supplier $i\in\mathcal{S}$ provides. We also define the prices of the products for the consumers as $\pi_j:=\pi_{p(j)}$ and we define $\pi_k:=\sum_{p\in \mathcal{P}}\gamma_{k,p}\pi_p$ as the technology price/value (this is the weighted sum of its input and output products, weighted by the transformation factors). If we now define the profit functions:
\begin{subequations}
\label{SPS}
\begin{align}
\phi_{i} &:= (\pi_{i}-\alpha_{i})s_i, \quad i \in \mathcal{S}
\label{SPT}\\
\phi_{k} &:= (\pi_{t}-\alpha_{k})t_k, \quad k \in \mathcal{T}
\label{SPD}\\
\phi_{j} &:= (\alpha_{j}-\pi_{j})d_j, \quad j \in \mathcal{D}, 
\end{align}
\end{subequations}
we can see that the Lagrange function can be rewritten as:
\begin{equation}
\begin{aligned}
\mathcal{L}(s,d,t,\pi) = \sum_{j \in \mathcal{D}}{\phi_{j}} + \sum_{i \in \mathcal{S}}{\phi_{i}} + \sum_{k \in \mathcal{T}}{\phi_{k}}. 
\end{aligned}
\end{equation}
Because the optimization problem is a linear program, we can assume that strong duality holds and thus its solution can be determined by solving the Lagrangian dual problem: 
\begin{equation} 
\label{MLD}
\begin{aligned}
\max_{\pi}\max_{(s,d,t) \in \mathcal{C}} \mathcal{L}(s,d,t,\pi)
\end{aligned}
\end{equation} 
Here, the set $\mathcal{C}$ captures the capacity constraints \eqref{MS},\eqref{MD}, and \eqref{MT}. From these observations, we can see that the optimization model is seeking to maximize the total profit of all system elements (suppliers, consumers, and technologies). Moreover, because strong duality holds, we have that the optimal allocations and prices are such that total surplus equals the total profit: 
\begin{align}\label{eq:totalsurplusprofit}
\sum_{j \in \mathcal{D}}{\alpha_{j}}d_{j} - \sum_{i \in \mathcal{S}}{\alpha_{i}}s_{i} - \sum_{k\in\mathcal{T}}{\alpha_{k}t_{k}}=\sum_{j \in \mathcal{D}}{\phi_{j}} + \sum_{i \in \mathcal{S}}{\phi_{i}} + \sum_{k \in \mathcal{T}}{\phi_{k}}.
\end{align}
The total surplus can thus be interpreted as the {\em system-wide profit}; this profit in turn can be interpreted as the {\em total wealth} created by the system (i.e., the system can be interpreted as a market or an economy). 
\\

Lagrangian duality analysis also reveals the role of the dual variables in remunerating elements. Specifically, the profit of the supplier is given by difference between the actual payment ($\pi_is_i$) and the expected cost ($\alpha_is_i$); as such, the supplier desires that $\pi_i\geq \alpha_i$ (this guarantees non-negative profit) and that $\pi_i$ is as large as possible (to maximize its profit). For the consumers, we have that the profit is given by the difference between the expected cost ($\alpha_jd_j$) and the actual payment ($\pi_jd_j$); as such, we the consumer desired that $\pi_j\leq \alpha_j$ (this guarantees non-negative profit) and that $\pi_j$ is as low as possible (to maximize its profit). For the technologies we have a more interesting case; specifically, for the profit to be non-negative, we require that $\pi_k\geq \alpha_k$; by definition, we have that $\pi_k=\sum_{p\in \mathcal{P}}\gamma_{k,p}\pi_p$. As such, it is important to observe that we need that the {\em value created by a technology} needs to be higher than its operating cost $\alpha_k$. This can only be achieved if the generated output products have a higher value than the consumed input products (weighted by the transformation factors). This insight is important, as it clearly shows how technologies are {\em wealth creators} (generate higher-value products from lower-value products). Moreover, any technology that cannot achieve wealth creation (due to mismatches in input/output costs or technologies inefficiencies) will simply not participate in the economy. 
\\

We note that the prices of different products obtained from the model are inherently linked; this enables to capture dependencies in product values (e.g., steam value is linked to the electricity value). This is an important property of the model and enables its use in other systems such as district utility systems, in which the price of thermal energy is linked to that of electrical energy \cite{sjodin2004calculating,li2015review}.
\\

Additional insight into the prices $\pi_p,\; p \in \mathcal{P}$ generated by the optimization problem can be determined by posing its dual problem:
\begin{subequations} 
\label{MDual}
\begin{align}
\min_{\pi,\lambda} \quad \sum_{i \in \mathcal{S}}{\bar{s}_{j}\lambda_{i}} + \sum_{j \in \mathcal{D}}{\bar{d}_{j}\lambda_{j}} + \sum_{k\in\mathcal{T}}{\bar{t}_{k}\lambda_{k}}
\label{MS_dual}\\
\text{s.t} \qquad \pi_{i} - \lambda_{i} \ge \alpha_{i}, \quad i \in \mathcal{S}, \quad (s_{i})
\label{MD_dual}\\
\qquad \qquad \pi_{j} - \lambda_{j} \le \alpha_{j}, \quad j \in \mathcal{D}, \quad (d_{j})
\label{MT_dual}\\
\qquad\qquad \pi_{k} - \lambda_{k} \ge \alpha_{k}, \quad k \in \mathcal{T}, \quad (t_{k})
\end{align}
\end{subequations} % show dual to get price bounds

This formulation reveals that the prices must be {\em bounded} by the offered values. Here, we have that $\lambda_i,i \in \mathcal{S}$, $\lambda_j, j \in \mathcal{D}$, and $\lambda_k, k \in \mathcal{T}$ are the dual variables assigned for the upper bound constraints of the suppliers, consumers, and technologies.  These dual variables must be positive; as such, they indicate that the supplier and technology prices must be higher than their offering values, while the prices for the consumers must be lower than their offering costs. This again indicates that prices generated by the model are such that all system elements generate a profit (or at least break even). This also indicates that the system-wide profit must be non-negative; specifically, the {\em entire} system must generate wealth (or at least break even); in other words, it makes no sense to have an economy that does not generate value. 
\\

Note also that an outcome of the optimization model is that all allocations are zero $(s,d,t)=(0,0,0)$. This can occur, for instance, if the supplied products and technology services are too expensive compared to values of requested products (and thus the system cannot generate wealth). As such, the system-wide profit can be seen as a measure of economic efficiency (a total profit of zero indicates the economy/system/value-chain is not economically efficient). Another possible outcome of the model is that the allocations of {\em some specific technology pathways} are zero; this would indicate that such pathways do not generate wealth. This also indicates that the model has the goal of identifying pathways that generate the most value (and eliminate the ones that do not generate value); as such, the model can be used to quickly identify promising pathways. From all these observations, we can see that the proposed model makes economic sense; this is important, as it facilitates analysis and interpretation of model outcomes. For instance, it is possible to use the framework to explore what values of technology efficiencies can help activate pathways and generate wealth. 
\\

\begin{figure}[htp!]
\center{\includegraphics[width=0.7\textwidth]{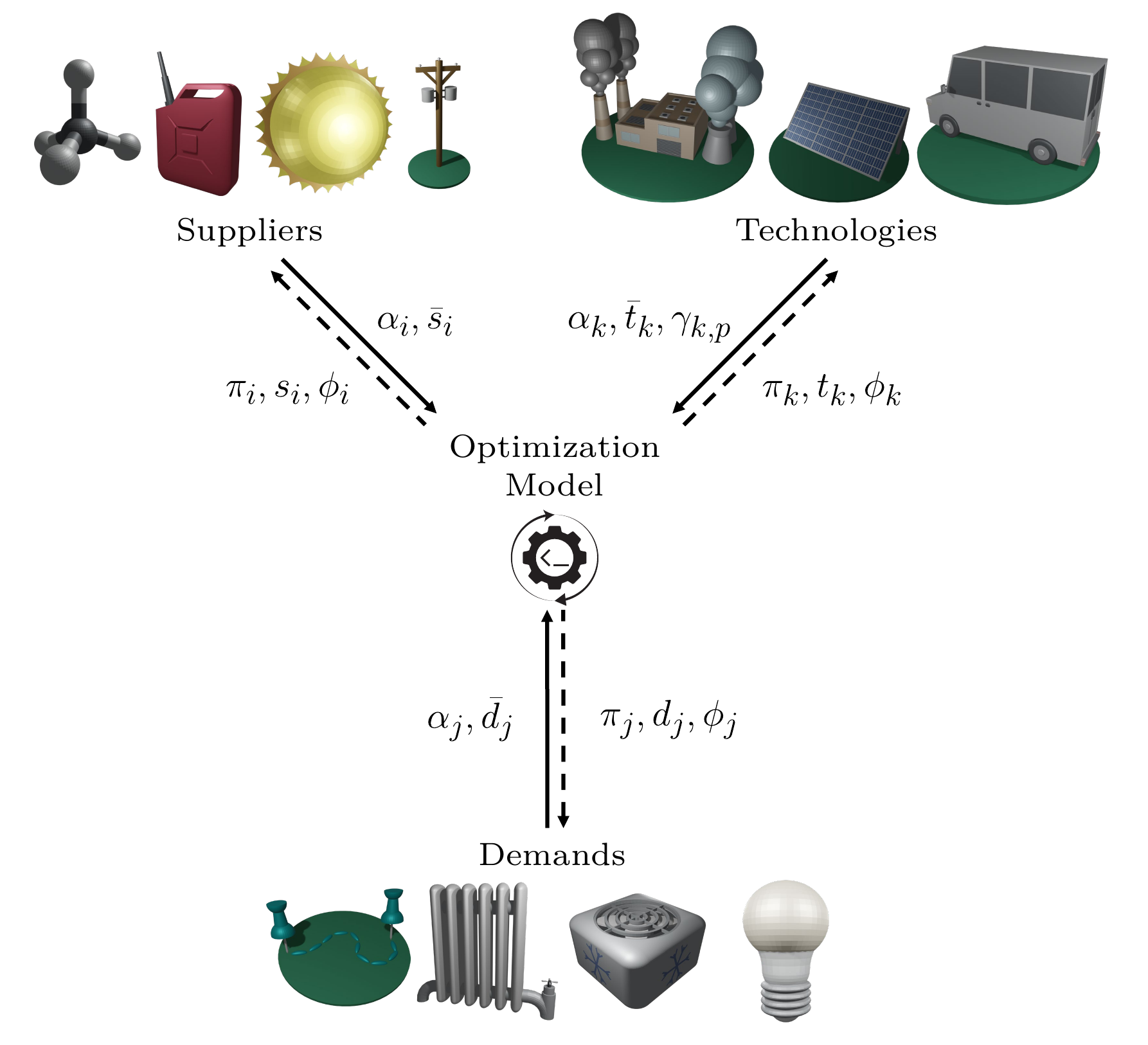}}
\caption{Schematic representation of optimization model, indicating elements involved as well as model inputs/outputs.}
\label{fig:model}
\end{figure}

The primal formulation \eqref{MPrimal} of the optimization provides the physical allocations of the products, while the dual variables $\pi_p,\, p \in \mathcal{P}$ from the dual formulation \eqref{MDual} provides economic allocations of the products. These formulations together make up the optimization model which takes in information from system elements, maximizes the total profit, and returns the physical and economic allocations to the elements as described by Figure \ref{fig:model}.
\\

The model proposed treats products $p\in\mathcal{P}$ as {\em abstract objects}; as such, these can represent any type of asset that is exchanged in the system (e.g., material, energy, labor, mileage, computing service). For instance, vehicles could take in fuel ($L~gas$) to produce distance ($km$) and emissions ($kg~CO_2$). The abstract product objects can also be used to represent a discrete spectrum of products (e.g., green, blue, and brown hydrogen or water of different types).
\\

The model also treats values for products as abstract objects and we can use these to capture a range of behaviors. For instance, offered values for suppliers ($\alpha_i$) and consumers ($\alpha_j$) are typically positive or non-negative. The price bounds obtained from duality indicate that, in such a case, the product prices will be positive. However, the offering values can can also be allowed to be negative; this provides flexibility to capture policy and taxes or negative impact of byproducts from the system (e.g., waste generation). For instance, one of the pathways shown in Figure \ref{fig:cars} produces a waste (CO\textsubscript{2}) as a byproduct of transportation. This waste CO\textsubscript{2} has to be taken by the environment; as such, we can interpret the environment as a consumer that is willing to take the waste product but at a negative cost (the negative cost can be interpreted as a carbon tax). This is analogous to how landfills and wastewater treatment facilities operate (they charge a tipping fee to take on waste). We can also envision a supplier that has a product that they want to get rid of; as such, they can offer the product at a negative cost (this is common when dealing with undesired waste). For instance, we are willing to pay a fee for garbage to be taken away from our homes (as garbage is an undesired product). Accounting for undesired outcomes of the economy (via negative prices) is important in capturing social or environmental impacts of the system. 
\\

We highlight that the optimization model does not enforces satisfaction of demands explicitly ($d_j=\bar{d}_j$); instead, the model will only decide to satisfy demands if this makes economic sense (i.e., if it maximizes the total profit). For instance, it might be possible that forcing the satisfaction of a demand of a non-valuable product will incur significant costs that render the system non-profitable. It is also worth highlighting that the model will select pathways between basic products and final products that make the most economic sense; as such, the model is fully driven by economic incentives. It is possible to modify the model to enforce demand satisfaction in a couple of ways. For instance, we can set a large offering value $\alpha_{j}^{\mathcal{D}}$ to demands that need to be satisfied; this is analogous to how electricity markets operate (a high value to electricity demands is set to ensure delivery). It is also possible to enforce strict constraints on satisfaction of demands, but such constraints can be shown to introduce economic inefficiencies \cite{tominac_economic_2021}, because this can force the selection of pathways that are non-profitable. 
\\

The proposed framework is useful as a screening tool because typical applications contain a large number of technology pathways. For instance, a superstructure of all possible pathways for a specific application is shown in Figure \ref{fig:SS}. This graph contains the current infrastructure and potential new pathways to reach decarbonization for a university campus (as we will discuss in our case study). For this example, electrolysis, fuel cells, electric vehicles, hydrogen vehicles, hydrogen generators, and combined heat and power systems can be added as potential technologies. This could become even more complex as new technologies are developed or new supplies are discovered or demands shift towards new products (thus expanding the boundary of the system). The bottom right demand is for CO\textsubscript{2} which has to be accounted for analysis, as this is an undesired byproduct of the economy/system. The proposed optimization model can be used to screen technology pathways that make most sense under different market conditions (e.g., under different carbon tax scenarios).
\\

\subsection{Investment Model}

It is often of interest to restrict the selection of technology pathways based on investment costs or budgets that might be available. We can easily extend the proposed model to account for this; we will use this formulation to study how available investment budgets can constraint certain technology pathways, as well as to understand how we should prioritize technologies. We define a binary variable $y_k \in \{0,1\},~k \in \mathcal{T}$ to represent if a technology is selected with an investment cost of $\beta_k\in \mathbb{R}$. This gives the formulation: 
\begin{subequations} 
\label{IPrimal}
\begin{align}
\label{IObj}
\max_{s,d,t,y} &\quad \sum_{j \in \mathcal{D}}{\alpha_{j}}d_{j} - \sum_{i \in \mathcal{S}}{\alpha_{i}}s_{i} - \sum_{k\in\mathcal{T}}{\alpha}_{k}t_{k}\\
\label{IBalance}
\text{s.t} &\quad \sum_{i \in \mathcal{S}_p}{s_i} + \sum_{k \in \mathcal{T}_{p}}{\gamma_{k,p}t_{k}} - \sum_{j \in \mathcal{D}_{p}}{d_{j}}=0, \quad p \in \mathcal{P}\\
\label{IS}
&\quad 0 \le s_{i} \le \bar{s}_{i}, \quad i \in \mathcal{S}\\
\label{ID}
&\quad 0 \le d_{j} \le \bar{d}_{j}, \quad j \in \mathcal{D}\\
\label{IT}
&\quad 0 \le t_{k} \le \bar{t}_{k}\cdot y_{k}, \quad k \in \mathcal{T}\\
\label{IB}
&\quad \sum_{k \in \mathcal{T}}{ \beta_{k}\cdot y_k} \le \bar{\beta}
\end{align}
\end{subequations}
The constraint \eqref{IB} enforces that the total investment cost is less than or equal to a given budget $\mathcal{\bar{\beta}}\in \mathbb{R_+}$. The maximum technology input constraint \eqref{IT} is augmented such that if the technology is not built ($y_k = 0$) then the corresponding input value will be zero ($t_k = 0$). For technologies that already in place, we simply set the corresponding binary values to one. In our modeling approach, we can capture economies of scale by defining technologies of different capacities. The budget constraint will enable studies to determine what budget is required to transition to different pathways. We note that it is also possible to incorporate the investment costs into the objective and let the model decide what investments are needed to maximize the total surplus (without any constraints in budget). In this context, it is important to ensure that operating costs and investment costs are on the same time basis (e.g., annualized costs). The investment model reduces to a management model if the binary variables are fixed. We can thus think of the investment model as a framework that can be used to study how strategic deployment of technologies can be used to activate an economy and generate wealth.
\\

%%%%%%%%%%%%%%%%%%%%%%%%%%%%%%%%%%%%%%%%%%

\section{Case Study: Decarbonizing a University Campus}

We consider a case study for a {\em simulated} university campus to illustrate the modeling capabilities of the framework and the types of insights that it can provide. The model was developed using basic data and assumptions that aim to capture general features of a prototypical, large campus (inspired by the UW-Madison campus); however, the model is not a real and validated representation of the actual system. As such, results and conclusions reached in this study are only intended to illustrate the capabilities of the model, and should not be used as actual recommendations. 
\\

\begin{figure}[!htp]
\center{\includegraphics[width=0.7\textwidth]{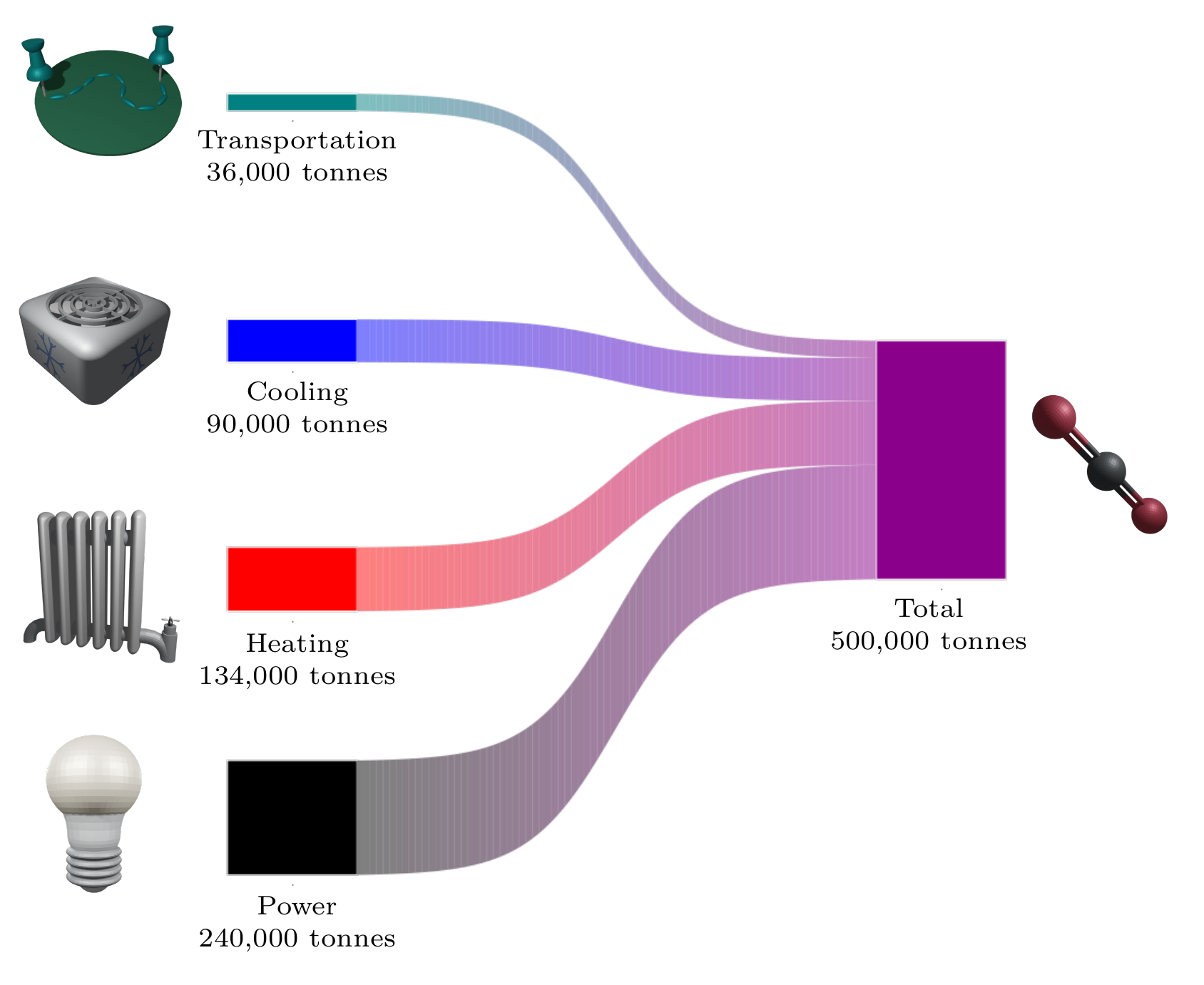}}
\caption{Sankey diagram summarizing the annual CO\textsubscript{2} emissions and their sources for a university campus. Generation of electrical power and heating are the dominant sources.}
\label{fig:Campus CO2}
\end{figure}

We consider a campus that serves 48,000 students and 24,000 staff members on a 3.6 km\textsuperscript{2} campus with 2.3 million m\textsuperscript{2} of buildings \cite{noauthor_university_nodate,frank_university_2022}. Under existing operations (and associated technology pathways), the university outputs 500,000 tonnes of CO\textsubscript{2} annually to meet demands of four key products (see Figure \ref{fig:Campus CO2}) \cite{frank_university_2022}. The carbon footprint is equivalent to the annual emissions of 34,000 average citizens in the US \cite{noauthor_world_2020} (thus comparable in magnitude to the campus population). To decrease the CO\textsubscript{2} footprint, we must identify new technology pathways that can satisfy the product demands (campus services) while reducing CO\textsubscript{2} emissions. 
\\

\subsection{Existing Technologies}

We assume that the campus purchases electricity from the local utility with a special rate scheme due to the large demand of the university of ~400,000 MWh per year \cite{frank_university_2022}. The university has some small renewable energy projects to generate electricity, but these total to only ~30 MWh per year. These are primarily rooftop solar photo-voltaic solar panels. The rest of the power comes from the power grid. The largest power sources for electricity in the state are coal and natural gas at 42\% and 34\%, respectively \cite{noauthor_eia_2022-1}. There is one large nuclear plant supplying 15\% of electricity and hydroelectric supplies 4\%; the remaining 5\% is composed of non-hydroelectric renewables. This leaves the grid electricity with an estimated CO\textsubscript{2} output of 560 kg CO\textsubscript{2}/MWh \cite{noauthor_eia_2022-1}. Grid electricity is represented by GP in Table \ref{tab:Acronyms} and in Table S4 and Table S5 (for detailed technology specifications), which takes in MWh grid and converts it into MWh and kg CO\textsubscript{2}. The grid is represented as a technology, to account for generated CO\textsubscript{2} emissions associated with the use of grid electricity. Operating costs for the technologies are estimated and are represented relative to the reference product.
\\

The campus has an annual average highest outdoor temperature of 35℃ and lowest temperature of -27℃, so both heating and cooling is required to keep buildings comfortable for teaching and research \cite{noauthor_noaa_2022}. To meet these heating and cooling demands, the campus uses utility network that delivers steam and chilled water to buildings \cite{jandl_lake_2019}. These services are supplied by a couple of heating and cooling plants and one combined heat and power (CHP) plant. These plants burn natural gas to produce steam, which is used to heat the buildings and are represented by SG\textsubscript{NG}. The plants also use refrigeration systems to chill water to cool the buildings, represented by WC. These plants, along with the CHP plant, meet the annual demand for the university of ~1 million tonnes of steam and ~24 million tonnes of chilled water to maintain buildings at a comfortable temperature as shown in Table S2.
\\

The natural gas-fired CHP plant has a rated electricity output of 150 MW, can produce 220 tonne/hr of steam, and 27,000 tonne/hr of chilled water \cite{noauthor_mge_nodate}. For this plant, the water chillers are powered by electricity and will be represented separately by WC. The steam and power output will be represented by SG\textsubscript{NG}.
\\

The campus has a fleet of 912 vehicles to maintain the campus and meet travel demands (measured as distance) \cite{frank_university_2022}. A total of 865 of the vehicles are exclusively gasoline powered with the rest being hybrid vehicles and just 2 fully electric vehicles. To compare different types of vehicles, there is a demand of distance (km) to meet demands rather than a demand for fuel for transportation purposes. These vehicles satisfy the 6 million km the university requires each year at a rate estimated by the university fleet rate \cite{noauthor_uw-madison_2020}. These vehicles are represented by GV with operating costs based on average values for internal combustion vehicles.
\\

\begin{figure}[htp!]
\center{\includegraphics[width=0.8\textwidth]{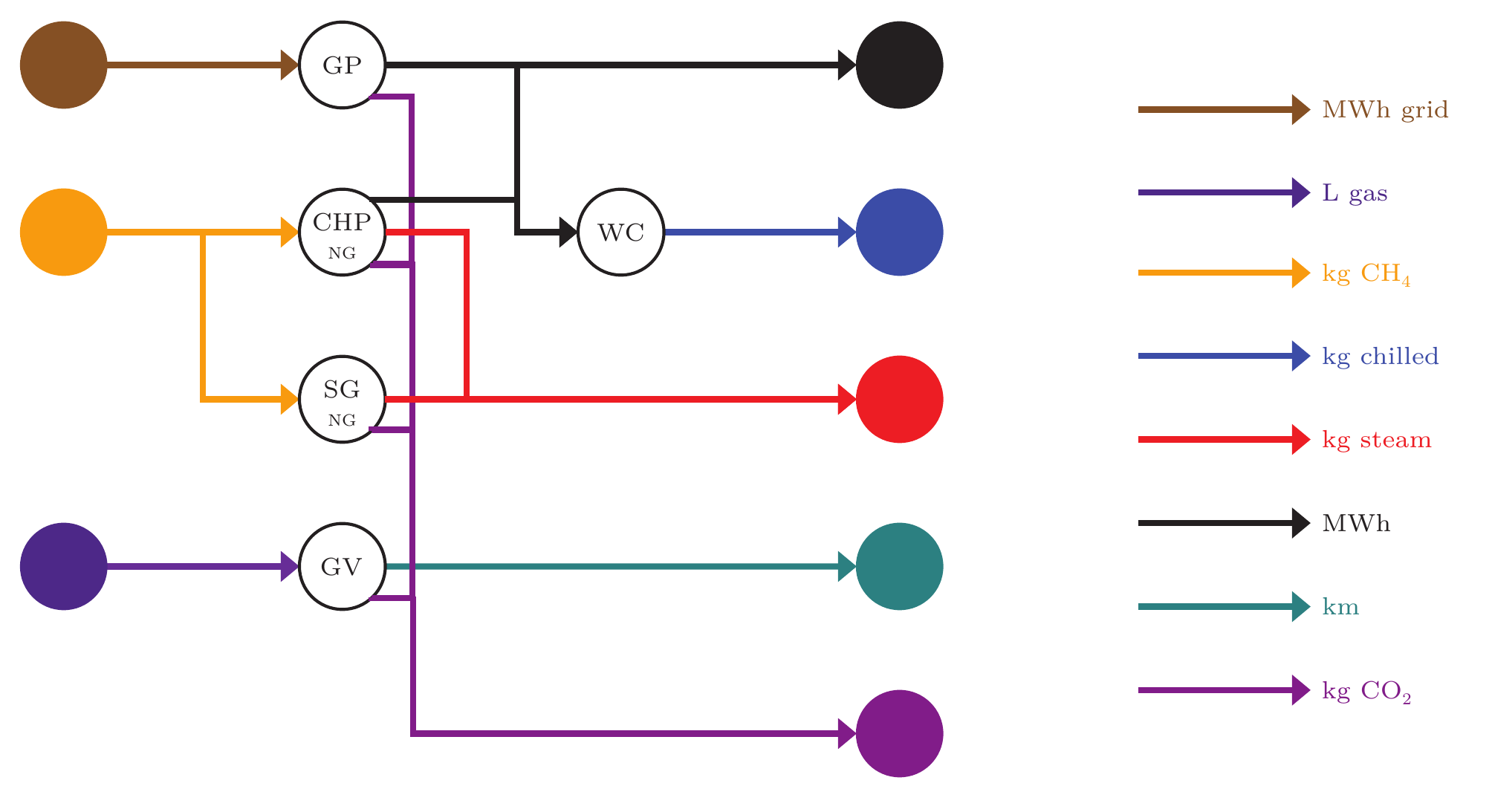}}
\caption{Superstructure graph showing currently-available suppliers, technologies, and demands for the university campus. See Table \ref{tab:Acronyms} for technology descriptions.}
\label{fig:Current}
\end{figure}

These technologies make up the current pathways that meet the campus demands and are displayed in Figure \ref{fig:Current}. Supplies of natural gas, gasoline, and grid electricity are added with pricing details found in Table S3. When implemented in the model, the quantities for these supplies are set to a sufficiently large value such that the technologies are not limited by the available supply. Full demand data can be found in S2. 
\\

While the case study campus occurs within one geographical location, this framework could be used to evaluate a network of campuses as well. For example, if the university had a satellite campus that can take power from the main campus but not steam and chilled water then there would be multiple types of steam and chilled water (e.g. kg steam and kg steam\textsubscript{satellite}). Similarly, if the power produced in the satellite campus could only be used within the satellite campus then the satellite campus would have a demand for MWh\textsubscript{satellite} and technologies that produce and use the power the in satellite campus would reflect this. To account for the power from the main campus being used in the satellite campus then there would need to be an additional technology that converts MWh to MWh\textsubscript{satellite}. 

\subsection{Potential New Technologies}

Electricity generation on campus produces the most CO\textsubscript{2}; to decarbonize generation, renewable electricity generation is required. To do this, the university would have to either wait for the power grid to use more renewable power or install their own renewable generators. However, renewable capacity depends on location and cannot be placed wherever is convenient. Further, the times that renewable power generation occurs does not always match the demand, meaning some energy would have to be stored for later use. Based on this constraint, the power produced from renewable power will be treated as a separate product ($MWh~renew$) that cannot be used directly to satisfy electricity demands. This can be thought of as a worst case scenario in which the renewable power generation ($MWh~renew$) never matches the electricity demand of the university and as such cannot be treated as electricity that can be used for any purpose ($MWh$). For the purposes of this case study, solar photovoltaic and wind power will be considered. Solar power will be represented by SP and wind power will be represented by WP to enable the use in the investment model and capture the operational costs of the two technologies using information from the U.S. Energy Information Administration \cite{noauthor_eia_2022}. These technologies are modular, meaning that additional capacity can be added in small increments for the same cost. As such, the investment cost and capacity reflect these aspects and, when implemented in the model, will have a set number of units available. These technologies take into consideration the required land usage, and an additional supply of land is added to account for this. This means that the technologies will be competing for a supply of land (in m\textsuperscript{2}). In other words, our model captures land as a resource that needs to be strategically used to various uses. This illustrates how the abstract products can be used to capture diverse assets.
\\

Hydrogen-based technologies have been proposed as an energy carrier with similar applications as natural gas. As such, hydrogen-based technologies could be critical to continue to meet the demands of the university. To be decarbonized, hydrogen can be produced through electrolysis using renewable power to create a CO\textsubscript{2} free energy carrier. As such, electrolysis will be represented by WS which takes in renewable power ($MWh~renew$) and produces hydrogen. The hydrogen could then be used by a fuel cell for electricity by FC, which is a modular technology as well \cite{marcinkoski_doe_2015, menezes_department_2020}. 
\\

The produced hydrogen can be used by the existing CHP plant by blending hydrogen and natural gas as represented by CHP\textsubscript{B}. This option has an investment cost of just 1\% of the initial investment cost of the turbines, because many natural gas powered turbines can handle hydrogen blends up to a volume fraction \cite{oberg_exploring_2022}. This study will use a blend of 20\% hydrogen by volume. Alternatively, the university could replace the existing turbines to run entirely off hydrogen with an associated investment cost of 25\% of the initial cost of the turbines of the CHP plant as represented by CHP\textsubscript{H\textsubscript{2}} \cite{oberg_exploring_2022}. Alternatively, hydrogen could be used by a fuel cell CHP plant which uses the high temperatures a fuel cell operates at to generate heat that could be used by the university as well as electricity. The fuel cell CHP will be represented by CHP\textsubscript{FC}. Hydrogen could also be used to decarbonize the heating demand of the university by implementing hydrogen steam generators as described by SG\textsubscript{H\textsubscript{2}}. This can be implemented without having to make any additional upgrades, unlike the CHP plant \cite{wang_efficiency_2022}.
\\

To reduce the CO\textsubscript{2} emissions of the university travel demand, alternative vehicles must be considered. Hybrid vehicles have 50\% higher fuel efficiency than traditional gasoline vehicles and will be represented by HV \cite{noauthor_nrel_nodate}. However, even the most fuel efficient gasoline vehicle still emits CO\textsubscript{2}. Electric vehicles could be a CO\textsubscript{2} free alternative if the electricity is sourced from renewable power such as wind and solar and will be represented by EV. The shift to electric vehicles has already been observed as more companies release new vehicles \cite{noauthor_us_nodate-2}. Hydrogen fuel-cell vehicles are another alternative, although they are less commercially available than electric vehicles. However, they have faster refueling times compared to the charging times of electric vehicles \cite{noauthor_us_nodate,noauthor_us_nodate-1}. The hydrogen vehicles in this study will be described by H\textsubscript{2}V. For all the alternative vehicle options, they will also be represented in a modular way meaning the purchase of one of these vehicles provides the capacity for 8000 km and the purchase of multiple vehicles will be required. 
\\

\begin{figure}[!htp]
\center{\includegraphics[width=1\textwidth]{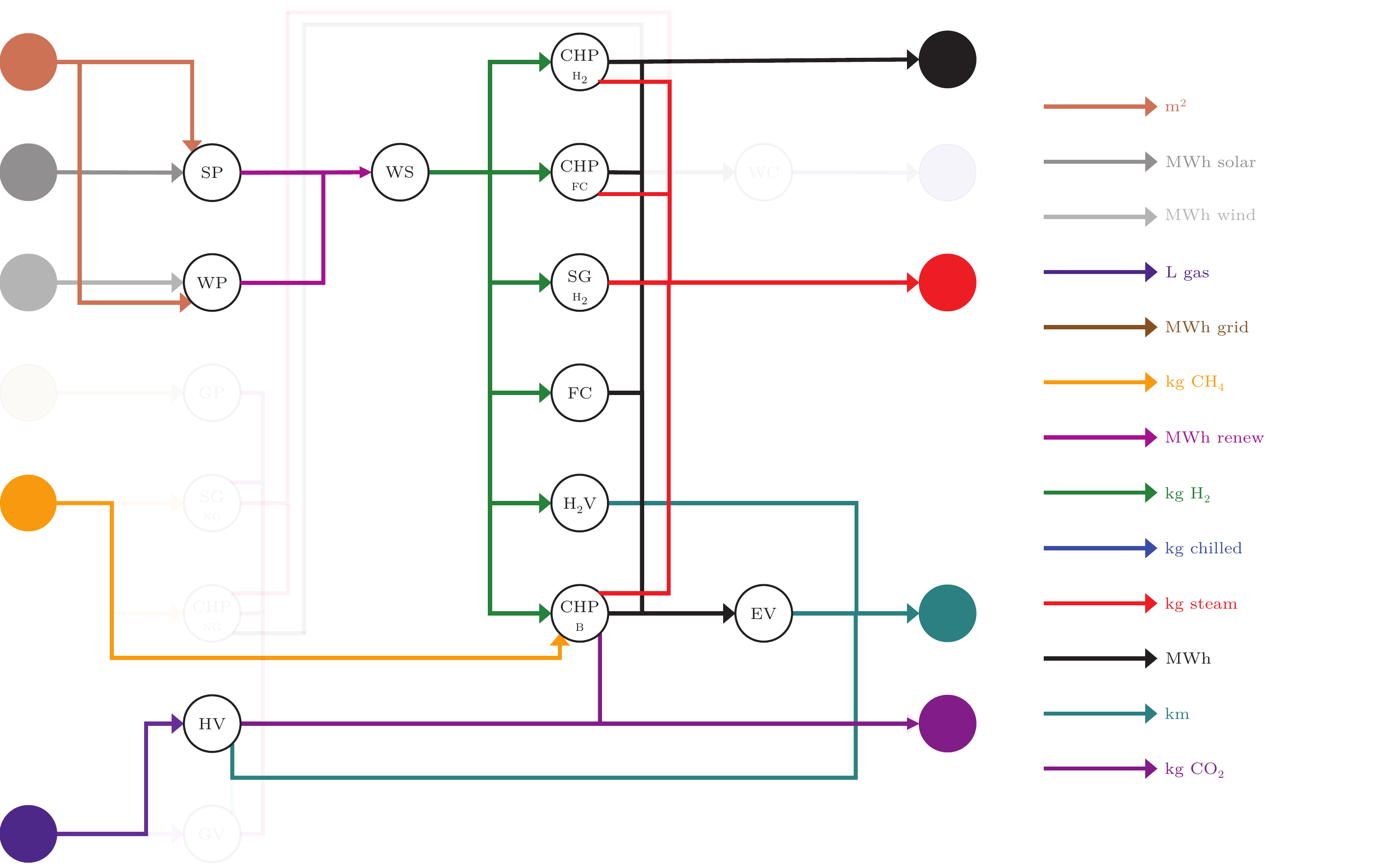}}
\caption{Superstructure graph of potential technology pathways that could be chosen by the university campus to decarbonize the system (current technologies have been blurred out to facilitate comparison).}
\label{fig:SS Potential}
\end{figure}

Modular technologies, such as vehicles and electrolysis, will have an associated number of units that are available in S7. While the single items like switching to a blend of hydrogen and natural gas feed for CHP are a single unit. This will allow the model to see what technology pathways are chosen at budgets that would not allow for a complete transition from, for example, gasoline vehicles to electric vehicles. All of the new technologies are shown in Figure \ref{fig:SS Potential}.
\\

All of the possible pathways from supply through technologies to demands associated with potential technologies are displayed in Figure \ref{fig:SS}. While all of the pathways are displayed here, a subset of them will be chosen by the model as being the optimal pathway. By varying some key attributes such as a carbon tax or the budget, we can determine what technologies and products will have the most impact under those conditions. 
\\

While the focus of the university campus presented here is on the addition of modular technologies and the upgrading of non modular technologies, non modular technologies could be added as potential technologies as well. This would require specifying the technology operating cost, investment cost, and conversion at multiple sizes. For example, if the university were to consider adding a natural gas fired power plant to make up for the additional electricity use due to the switch to electric vehicles, then the natural gas fired power plant would have to be specified at multiple sizes and added as potential technologies. For expansion of non modular technologies this would occur similarly with specification of the expansion at different sizes. This approach does require more knowledge into the scaling costs of technologies compared to a modular technology which does not have scaling, but it can be done eve though it may be more difficult for emerging technologies that have not had as much study.

\begin{table}[!htb]
\caption{Acronyms for technologies considered in university campus study.}
\begin{center}
\begin{tabular}{|c c|} 
 \hline
 \textbf{Acronym} & \textbf{Meaning} \\ [0.5ex] 
 \hline
 CHP\textsubscript{NG} & Combined Heat and Power Plant (using natural gas)\\ 
 \hline
 SG\textsubscript{NG} & Steam Generator (using natural gas)\\ 
 \hline
 WC & Water Chiller\\ 
 \hline
 GP & Grid Power\\ 
 \hline
 GV & Gasoline Vehicle \\  
 \hline
 SP & Solar Power (Photovoltaic) \\ 
 \hline
 WP & Wind Power \\ 
 \hline
 WS & Water Splitting (Electrolysis) \\ 
 \hline
 CHP\textsubscript{B} & Combined Heat and Power Plant (using blend of natural gas and hydrogen)  \\ 
 \hline
 CHP\textsubscript{H\textsubscript{2}} &  Combined Heat and Power Plant (using hydrogen) \\ 
 \hline
 CHP\textsubscript{FC} & Combined Heat and Power Plant (using Hydrogen Fuel Cell) \\ 
 \hline
 SG\textsubscript{H\textsubscript{2}} & Steam Generator (using Hydrogen) \\ 
 \hline
 FC & Hydrogen Fuel Cell \\ 
 \hline
 H\textsubscript{2}V & Hydrogen Fuel Cell Vehicle\\ 
 \hline
 EV & Electric Vehicle \\ 
 \hline
 HV & Hybrid Gasoline Vehicle \\ 
 \hline
\end{tabular}
\label{tab:Acronyms}
\end{center}
\end{table}

\subsection{Results}

\subsubsection{Management Model}

The purpose of the management model is to determine optimal technology pathways for the campus to meet its needs, without any restrictions on the budget needed to build such technologies. These optimal technologies will vary based on external conditions of the system (e.g., policy, markets, technology efficiencies). For instance, a key external condition is the cost of CO\textsubscript{2}. There is work to determine the externality of releasing CO\textsubscript{2} that places a cost of 113 USD per tonne CO\textsubscript{2}, although estimates vary \cite{wang_estimates_2019}. However, our work will aim to ask the question: What CO\textsubscript{2} cost/tax would incentivize the university to reduce or stop emitting CO\textsubscript{2}? The tax does not necessarily have to be interpreted as an externality (e.g., government-imposed), but can also be interpreted as an implicit value that the university is placing to its carbon footprint (which can potentially trigger negative public perception). 
\\

The carbon emissions are modeled as a demand for CO\textsubscript{2} that is charged at a negative bid value (it is a waste product). This can be interpreted in different ways; for instance, this can be interpreted as a tax that the government is introducing for any emissions generator by the system. Alternatively, this can be interpreted as a ``tipping fee" that the environment is charging the system to take its CO\textsubscript{2} waste. The CO\textsubscript{2} cost can also be interpreted as an internal (inherent) value that campus is placing on emissions; this is analogous to how companies that are currently seeking to decarbonize think about this waste (e.g., they are internally placing a negative value to CO\textsubscript{2} that might be implicitly connected to branding or public perception). This illustrates how the proposed modeling framework can help capture diverse scenarios of interest.
\\

While a negative bid value for CO\textsubscript{2} is used to encourage decarbonization, the other demands must be met. These inflexible demands for power, steam, chilled water, and distance must be met for the university to operate. As such, the demand bids will be set to an artificially large positive value. This means the university is willing to pay any cost to meet the demands. We can then rely the dual formulation of the problem \eqref{MDual} to discuss product values. This approach is similar to what is done with power system operation.
\\

\begin{figure}[!htp]
\center{\includegraphics[width=1\textwidth]{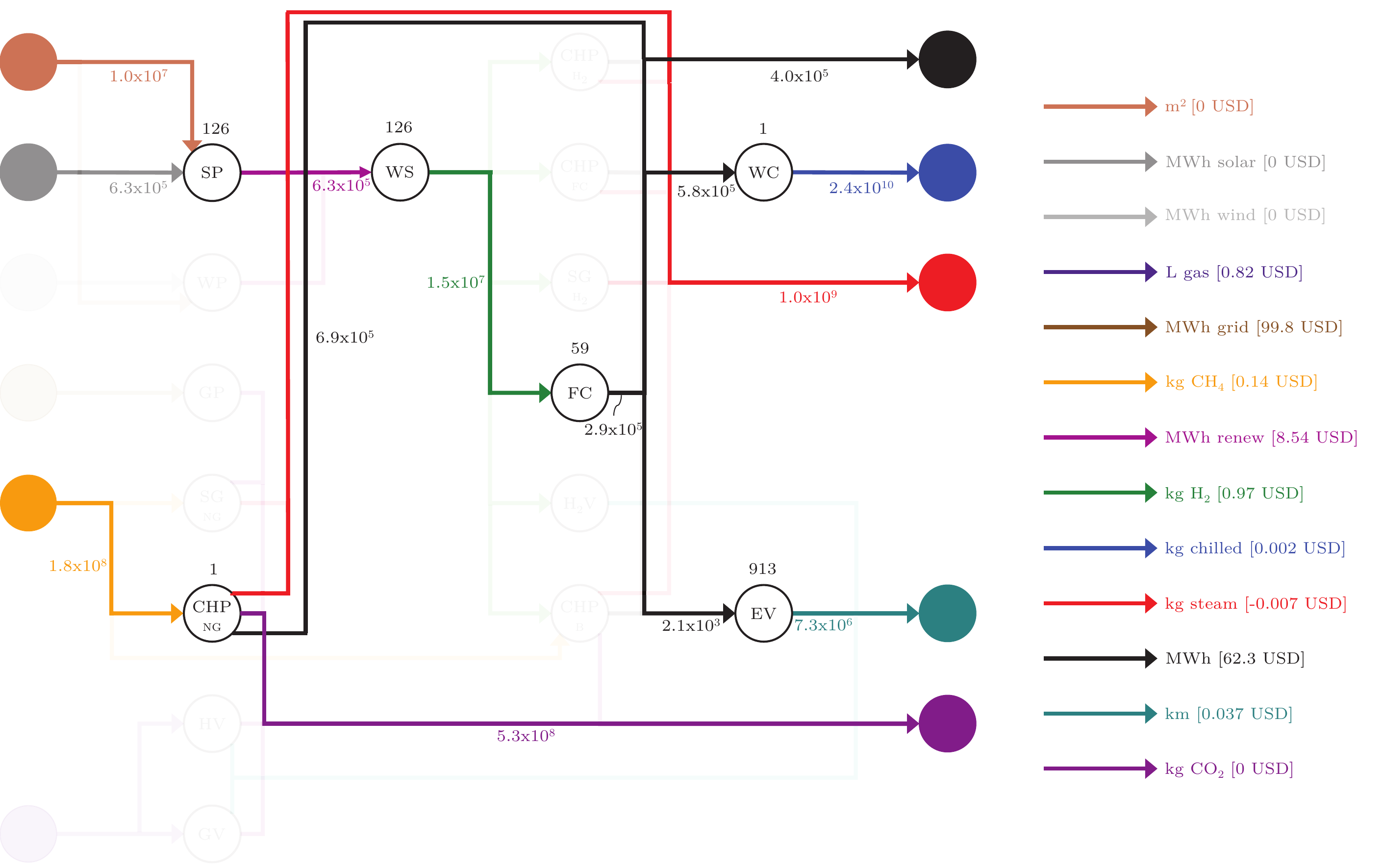}}
\caption{Optimal technology pathways under no CO\textsubscript{2} tax showing the technologies used satisfy the demands of university campus; pathways that are not selected are blurred out. Product prices are displayed in the legend, flows going in and/or out of each stakeholder are reported in the units of the legend, and the number of technology units required are reported above each technology.}
\label{fig:M1}
\end{figure}

To establish a benchmark, the pathway model with all potential and current technologies will first be optimized with no carbon tax. This leads to the technology pathway as shown in Figure \ref{fig:M1}. This pathway, however, still emits CO\textsubscript{2} because it is using the natural gas-powered CHP plant as much as possible. Meaning that without an incentive, the university has no need to fully decarbonize. The pathway, however, makes use of solar power to produce hydrogen to be used to supplement the CHP electricity for use in electric vehicles, to operate the water chiller, and meet the electricity demands of the university.
\\

\begin{figure}[!htp]
\center{\includegraphics[width=1\textwidth]{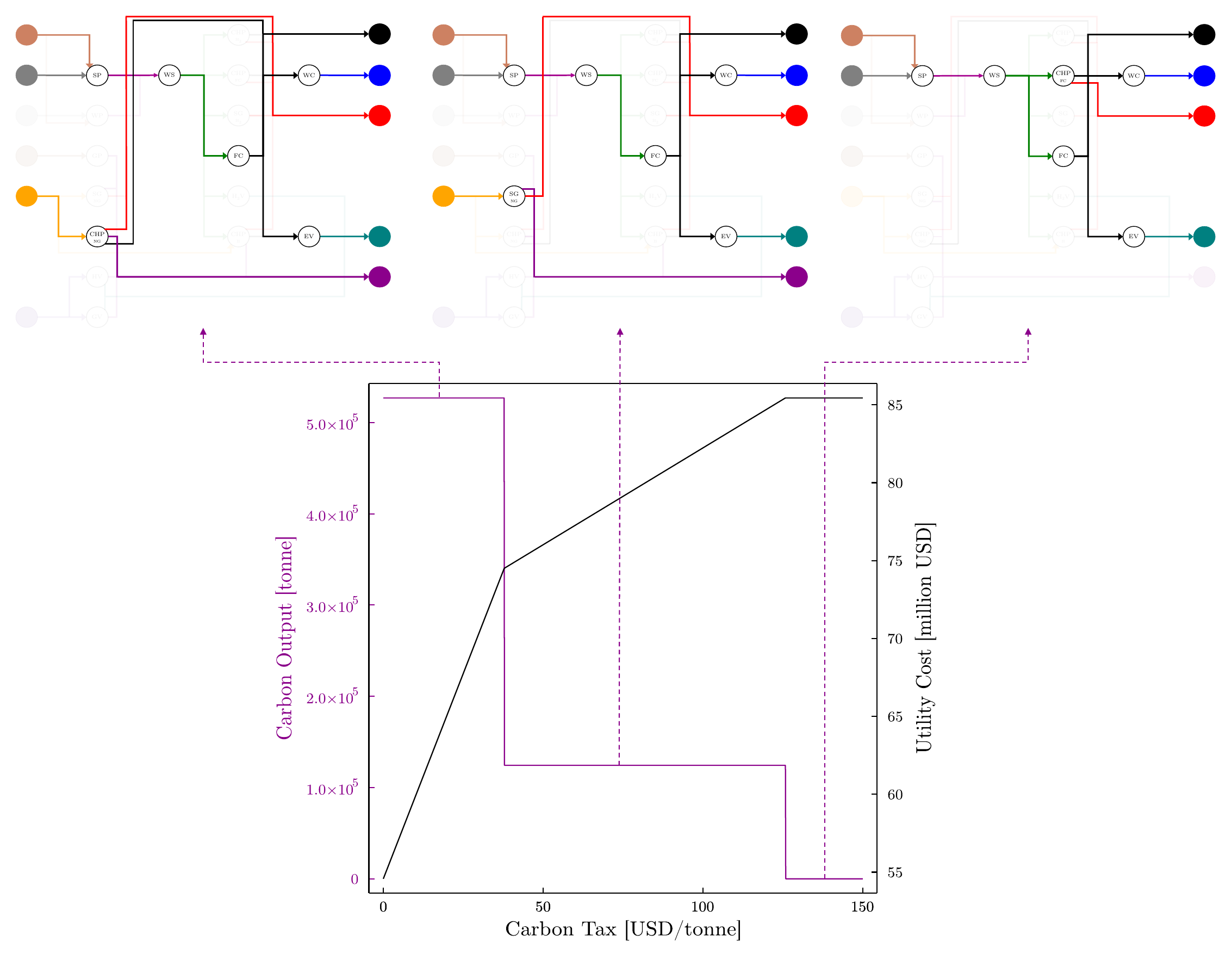}}
\caption{Dependence of CO\textsubscript{2} output and utility cost of the university on the CO\textsubscript{2} tax (when all possible technologies can be used). Results indicate that there are three different pathways that are activated at different tax levels. A version of this figure showing carbon output and utility costs in relative terms is displayed in Figure S1.}
\label{fig:MResults}
\end{figure}

To study the impact of the CO\textsubscript{2} tax ($\alpha_{kg~CO_2}$), it will involve discretizing a range of values and solving the management model for each value. This is computationally cheap as the optimization problem is a linear program and therefore quick to solve. The demand bid for CO\textsubscript{2} will be negative because it is being modeled as a tax; meaning the university will have to pay the environment/government to take the CO\textsubscript{2} it emits. The impact of this tax on the utility cost of the university was also determined; utility cost is the cost to operate the technologies, cost of supplies, and disposal cost of waste products (CO\textsubscript{2}), and is given by:
\begin{equation} 
\label{UBill}
\begin{aligned}
\sum_{k\in\mathcal{T}}{\alpha}_{k}t_{k} + \sum_{i \in \mathcal{S}}{\alpha_{i}}s_{i} - \alpha_{kg~CO_2}d_{kg~CO_2}
\end{aligned}
\end{equation}
Figure \ref{fig:MResults} shows the impact of varying the CO\textsubscript{2} tax on the CO\textsubscript{2} emissions and utility cost of the university. This shows three different technology pathways depending on the CO\textsubscript{2} tax. Below 45 USD per tonne, the pathway chosen is the same as in Figure \ref{fig:M1}. While the CO\textsubscript{2} output stays constant for each pathway, the utility cost increases when there is still CO\textsubscript{2} emissions because the university has to pay for the environment to take the emissions. 
\\

\begin{figure}[htp!]
\center{\includegraphics[width=1\textwidth]{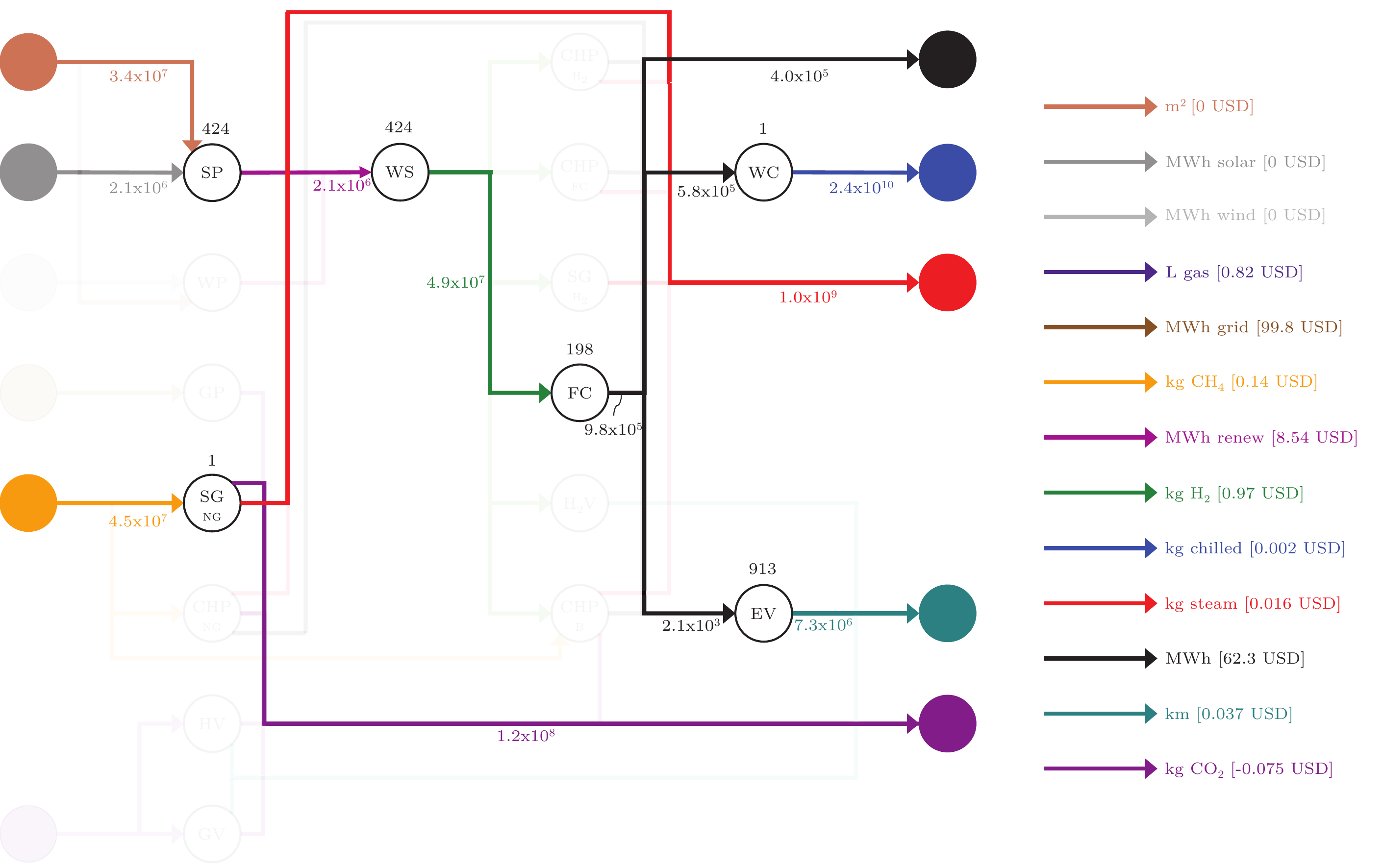}}
\caption{Optimal technology pathways for a CO\textsubscript{2} tax between 45 USD and 130 USD per tonne showing the technologies used satisfy the demands of university campus; pathways that are not selected are blurred out. Product prices are displayed in the legend for a CO\textsubscript{2} 0.075 USD per tonne, flows going in and/or out of each stakeholder are reported in the units of the legend, and the number of technology units required are reported above each technology. This pathway has decarbonized all of the university's demands besides the one for steam which still relies on natural gas fired steam generation.}
\label{fig:M2}
\end{figure}

Between ~45 USD and 130 USD per tonne is the pathway shown in Figure \ref{fig:M2}. This pathway is the one that would be chosen if the estimated externality of CO\textsubscript{2} was implemented as a tax. The pathway has a 77\% lower CO\textsubscript{2} output than in Figure \ref{fig:M1} but also has a higher utility cost. This pathway has fully decarbonized the demand for electricity, cooling, and transportation but still relies on natural gas for steam. 
\\

\begin{figure}[!htp]
\center{\includegraphics[width=1\textwidth]{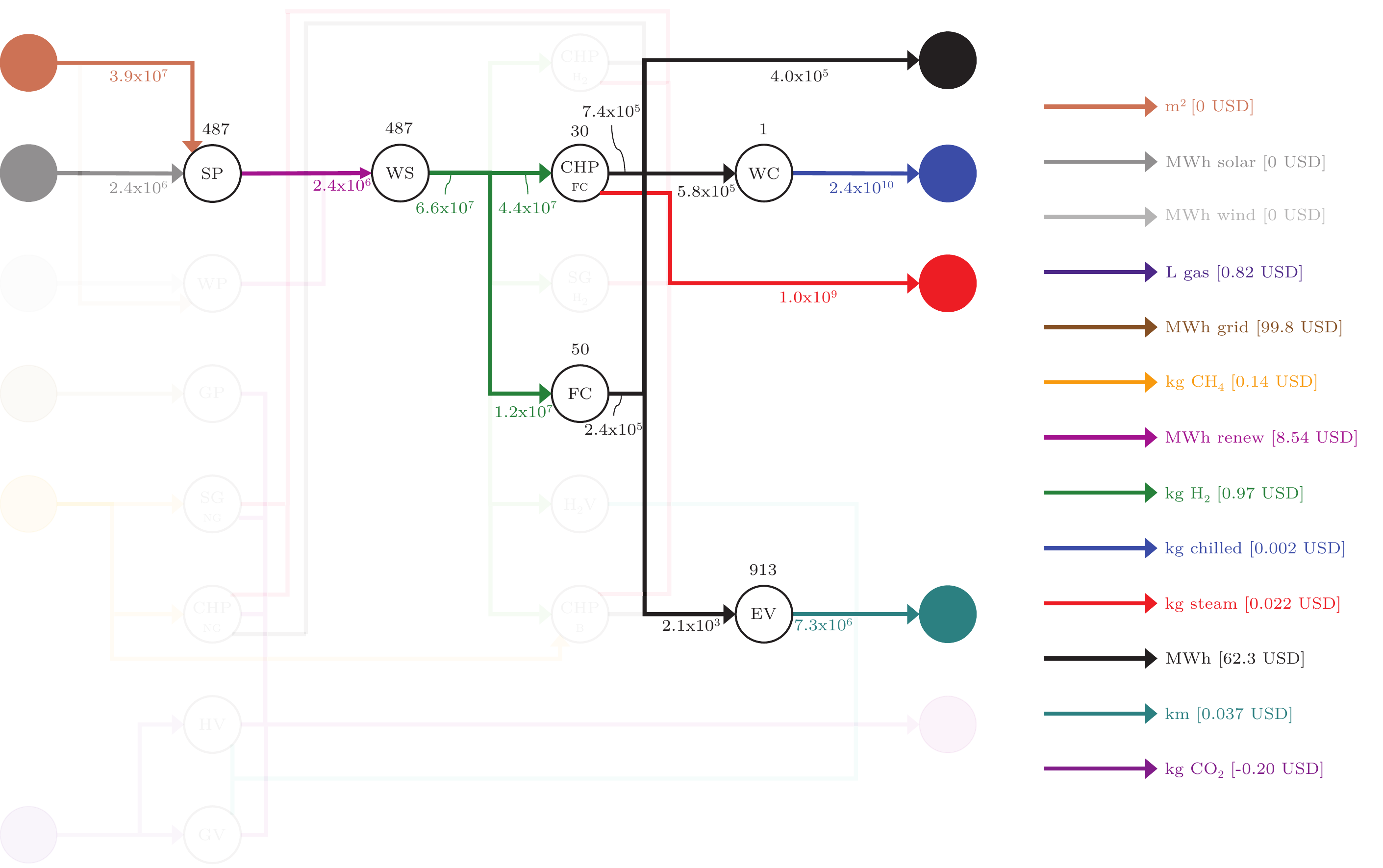}}
\caption{Optimal technology pathways for a CO\textsubscript{2} taxes above 130 USD per tonne showing the technologies used satisfy the demands of university campus; pathways that are not selected are blurred out. Product prices are displayed in the legend, flows going in and/or out of each stakeholder are reported in the units of the legend, and the number of technology units required are reported above each technology. This pathway has decarbonized all of the university's demands through the use of solar power and hydrogen fuel cell CHP and hydrogen fuel cells.}
\label{fig:M3}
\end{figure}

For CO\textsubscript{2} taxes above 130 USD per tonne, the optimal pathway is shown in Figure \ref{fig:M3}. This pathway has no CO\textsubscript{2} emissions because it has switched to using hydrogen fuel cell CHP to meet the demands for steam and some of the electricity. Fuel cells are used to supplement the electricity for use by the electric vehicles for distance, the water chiller for cooling, and the universities demand for electricity. This has the highest utility cost at ~55\% higher than when there is no CO\textsubscript{2} tax. This pathway also has the highest steam price as the fuel cell CHP plant is the most expensive to operate of the options for CHP, while all other prices remained constant between the three different pathways. 
\\

Switching from one of these technology pathways to another is not a simple task as it would involve investments and building. Investing any amount of money at a university can be a highly complex decision involving bids, board evaluation, state evaluation, fund procurement, and more. To advise a university seeking to decarbonize, the operating costs and investment amount for each pathway needs to be considered. The decision could be made, for example, by comparing the net present value of switching from the pathway in Figure \ref{fig:M2} to the pathway in Figure \ref{fig:M3} at a given CO\textsubscript{2} tax.
\\

The management model results identify that meeting the demand for steam to heat the buildings seems to be a major factor that drives technology selection. Additionally, the choice of using electric vehicles was consistent for the three different pathways because they have the lowest cost to operate per distance and the low cost of power from renewables. Solar power is chosen over wind power for three pathways because it has a lower operating cost per MWh produced than wind. With the use of solar power in each pathway, the hydrogen price remained constant at 0.97 USD/kg H\textsubscript{2} for the three different pathways. The electricity price also remained constant at 62.3 USD/MWh as it was set by the solar power to electrolysis to fuel cell pathway operating cost. This made the price for distance (km) and chilled water to remain constant. 
\\

Steam price did not remain constant across the three different pathways. Under no CO\textsubscript{2} tax (Figure \ref{fig:M1}) the steam price is -0.007 USD/kg. The negative value means that the natural gas fired CHP would be willing to pay for the steam to be taken because it already profits from electricity price being set by the fuel cell pathway. When the steam is instead produced by the natural gas fired steam generator (Figure \ref{fig:M2}) the steam price is no longer tied to the electricity price. The operating cost of the steam generator and the tax on the CO\textsubscript{2} released must be covered by the steam demand and a positive price of 0.016 USD/kg steam is observed at a CO\textsubscript{2} tax of 0.075 USD/kg CO\textsubscript{2}. When the steam is produced by the hydrogen fuel cell CHP and there are no CO\textsubscript{2} emissions (Figure \ref{fig:M2}) the steam price is 0.022 USD/kg steam. For this pathway, the positive price is a result of the revenue required from the steam to offset the higher operating cost of the hydrogen fuel cell CHP than the hydrogen fuel cell.
\\ 

The university modeled here benefits from economies of scale compared to a smaller university campus such as at a liberal arts college. This means that if a smaller university was modeled with different power procurement strategies such as purchasing power from the grid, installing a natural gas fired power plant, and installing solar power, then the conditions leading to the section of these technologies would be different than in the larger university presented in this case study. This is due to the different operating costs of the natural gas fired power plants per MWh between the two university campuses due to economies of scale. We would expect to see the transition to solar power from the natural gas fired power plant at a lower CO\textsubscript{2} value than for the larger university or even no use of a natural gas fired power plant if the operating costs were larger than for solar power.
\\

\subsubsection{Investment Model}

The investment model allows for the further consideration of a budget on choosing the technology pathways under different conditions. The key condition to incentives reducing CO\textsubscript{2} emissions will be a CO\textsubscript{2} tax, and analysis will involve solving the model for each value in a range. However, the impact of the budget has to also be considered. This means for each CO\textsubscript{2} tax, the investment model will have to optimize for each budget at each CO\textsubscript{2} tax. This is more computationally expensive than the management model because of the binary variables, but it is still readily solvable for each CO\textsubscript{2} tax and budget combination. 
\\

\begin{figure}[!htp]
\center{\includegraphics[width=0.7\textwidth]{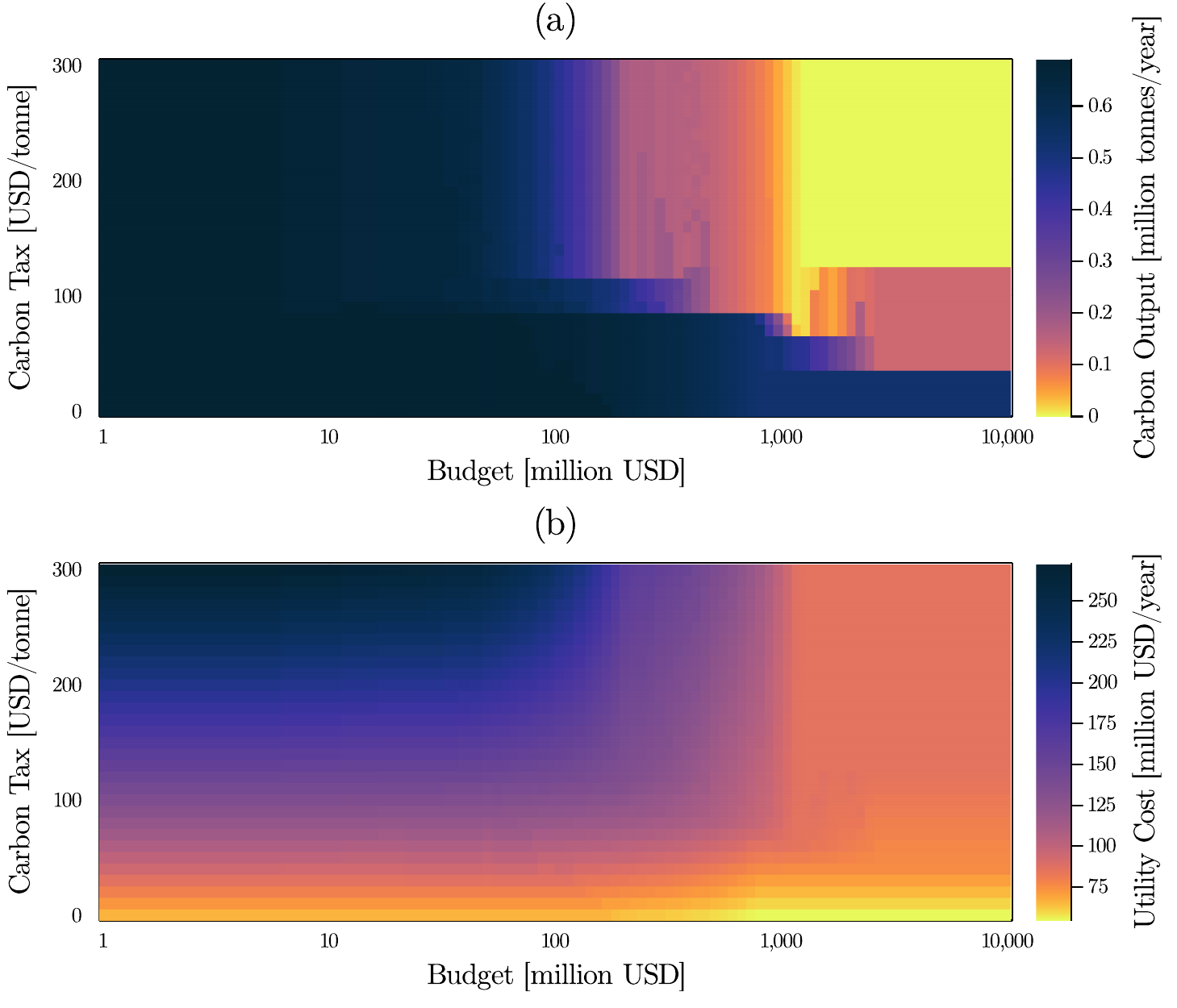}}
\caption{Dependence of the annual (a) CO\textsubscript{2} output and (b) utility cost of the university on the CO\textsubscript{2} tax and budget when all possible technologies are available for purchase and use. A version of this figure showing carbon output and utility costs in relative terms is displayed in Figure S2.}
\label{fig:IResults}
\end{figure}

Figure \ref{fig:IResults} shows the CO\textsubscript{2} emissions (a) and utility cost of the university (b) at varying CO\textsubscript{2} taxes and at varying investment budgets. When the budget is low, little can be done to shift away from the existing infrastructure so as the CO\textsubscript{2} tax increases, the emissions stay the same and the utility cost increases. This can be observed for all CO\textsubscript{2} taxes when the budget is less than ~80 million USD. While at sufficiently large budgets ($\ge$2.7x10\textsuperscript{9} USD), the potential pathways match the results of the management model and the same dependence on the CO\textsubscript{2} tax is observed. This occurs because the investment model has a high enough budget that it can choose any technology it desires, essentially eliminating the budget constraint. 
\\

\begin{figure}[!htp]
\center{\includegraphics[width=1\textwidth]{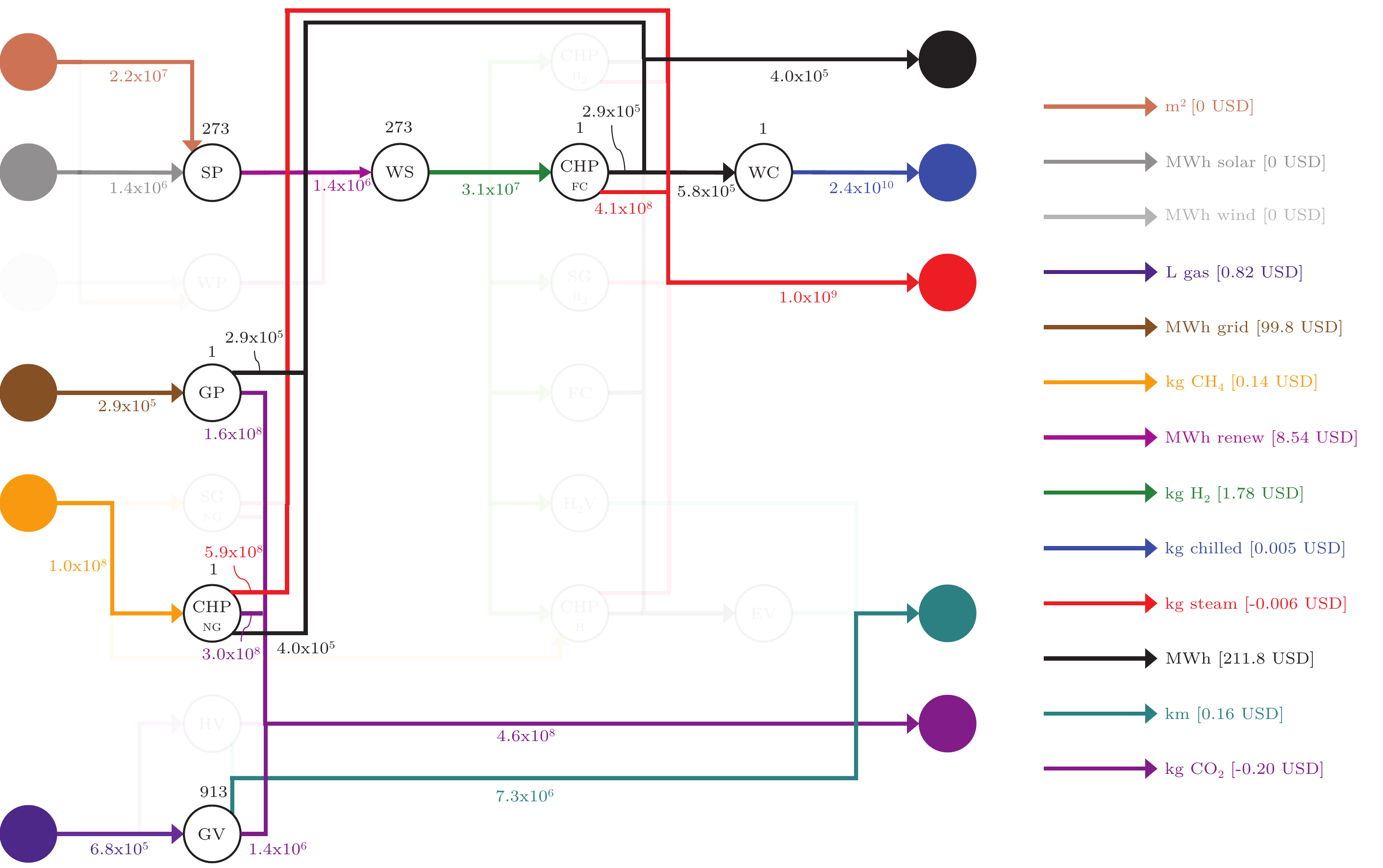}}
\caption{Investment model results at a CO\textsubscript{2} tax of 200 USD per tonne and a budget of 100 million USD showing the technologies used satisfy the demands of university campus; pathways that are not selected are blurred out. Product prices are displayed in the legend, flows going in and/or out of each stakeholder are reported in the units of the legend, and the number of technology units required are reported above each technology. This pathway makes use of some solar power to hydrogen for a fuel cell CHP unit however not demand is completely decarbonized.}
\label{fig:I1}
\end{figure}

As the budget increases above ~80 million USD, technologies are chosen to reduce the impact of the CO\textsubscript{2} tax. Figure \ref{fig:I1}, shows the pathway at a budget of 100 million USD and a CO\textsubscript{2} of 200 USD per tonne. This pathway is at a point where the budget is insufficient to allow for a complete switch to technology alternatives. Based on the management results, the optimal pathway at this CO\textsubscript{2} tax would be the one observed in Figure \ref{fig:M3}; however, the budget is not high enough for this pathway to be possible. Instead, the pathway chosen uses the existing gasoline vehicles to meet the distance demand, the natural gas powered CHP, and grid power which all produce CO\textsubscript{2}. There is some use of solar power for production of hydrogen which is then used by a single fuel cell CHP plant. The number of solar power units here are less than in the management model. The inability to offset the CO\textsubscript{2} tax results in a 2.4 times larger price for the power in this pathway than what the price would be without any budget constraint.
\\

Another point where the budget limits the pathway options is at a budget of 1x10\textsuperscript{9} USD and a CO\textsubscript{2} tax of 75 USD per tonne. At this CO\textsubscript{2} tax, the optimal pathway without any budget constraint would be the one displayed in Figure \ref{fig:M2} which had decarbonized everything besides the demand for steam. However, the budget limits the pathway to continue using the natural gas CHP plant and the gasoline vehicles. There is more use of solar power than in Figure \ref{fig:I1} due to the higher budget and more fuel cell CHP units, but it is still insufficient to meet the full demand for steam and electricity. Fuel cell units are used in this pathway to supplement the electricity from the natural gas and fuel cell CHP as well. This pathway has higher prices than in the budget unconstrained case with a 35\% higher power price. This is a result of the electricity price being impacted by the operating cost of the natural gas fired CHP and the tax on the CO\textsubscript{2} it emits. 
\\

\begin{figure}[htp!]
\center{\includegraphics[width=1\textwidth]{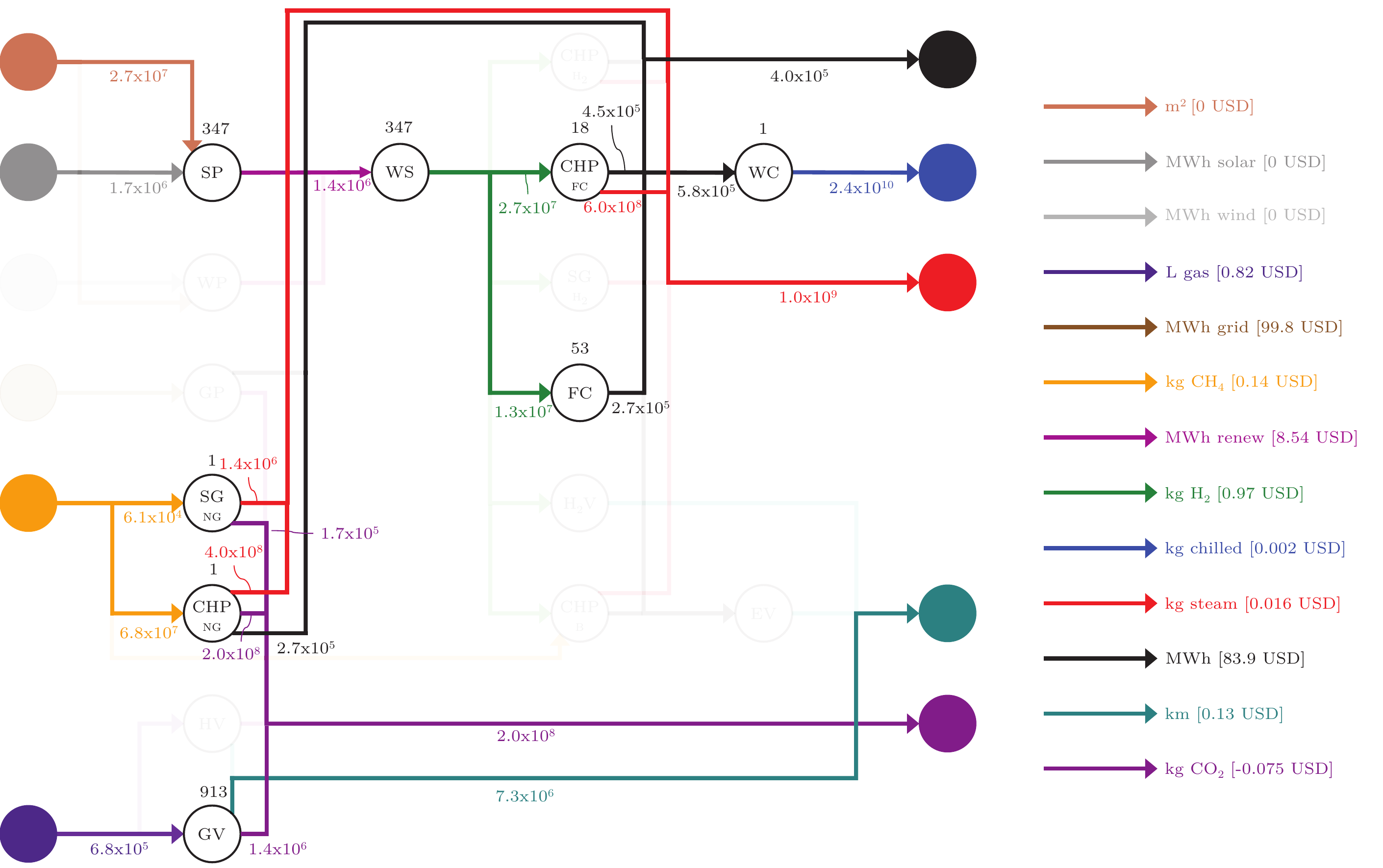}}
\caption{Investment model results at a CO\textsubscript{2} tax of 75 USD per tonne and a budget of 1x10\textsuperscript{9} USD showing the technologies used satisfy the demands of university campus; pathways that are not selected are blurred out. Product prices are displayed in the legend, flows going in and/or out of each stakeholder are reported in the units of the legend, and the number of technology units required are reported above each technology. This pathway uses no grid power but still uses natural gas fired CHP and steam generators. Solar power is used to produce hydrogen to be used by the fuel cell CHP and fuel cells.}
\label{fig:I2}
\end{figure}

\begin{figure}[htp!]
\center{\includegraphics[width=0.7\textwidth]{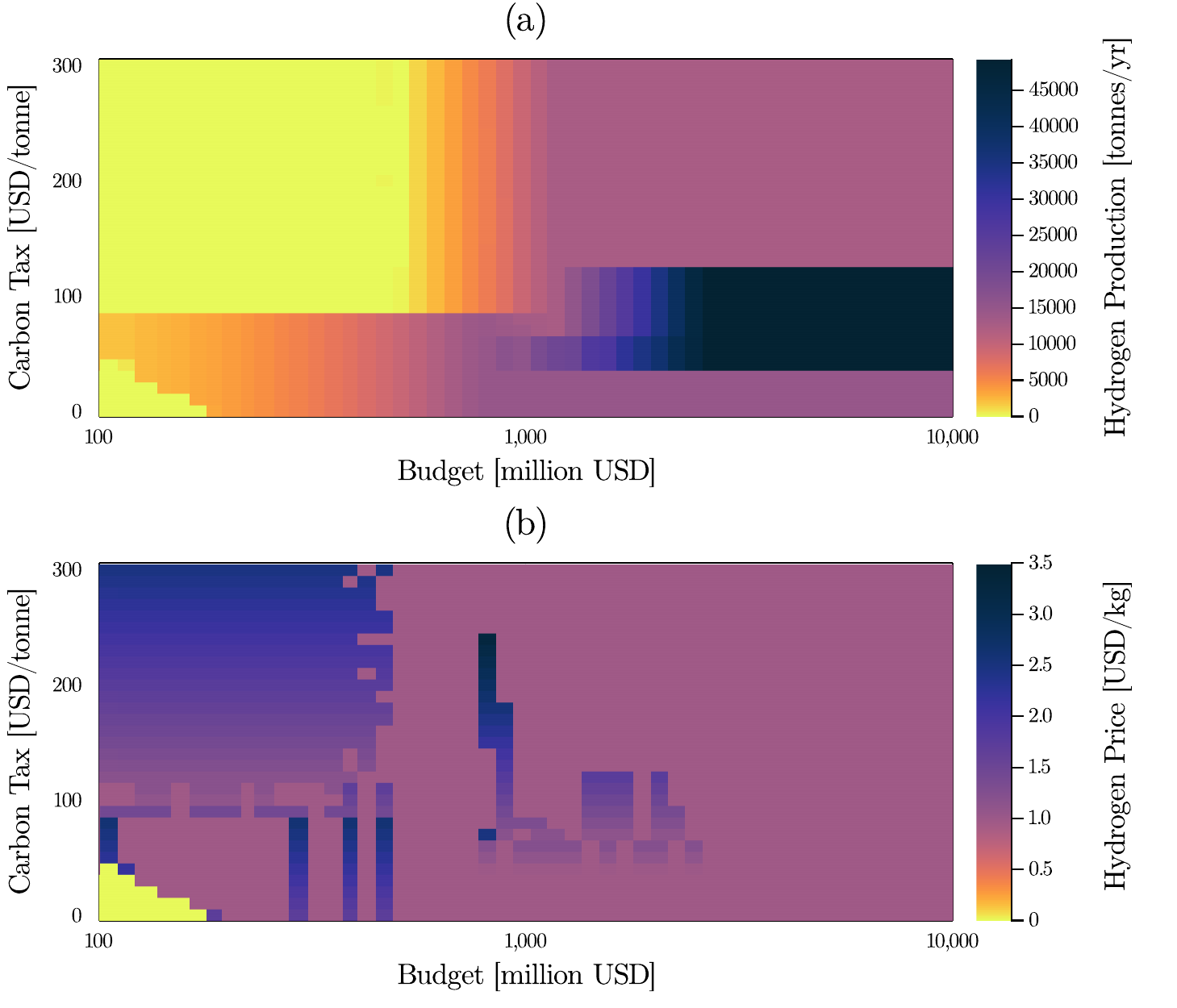}}
\caption{Dependence of (a) hydrogen production and (b) hydrogen prices on the CO\textsubscript{2} tax and budget when all possible technologies are available for purchase and use.}
\label{fig:IRP}
\end{figure}

\begin{figure}[htp!]
\center{\includegraphics[width=1\textwidth]{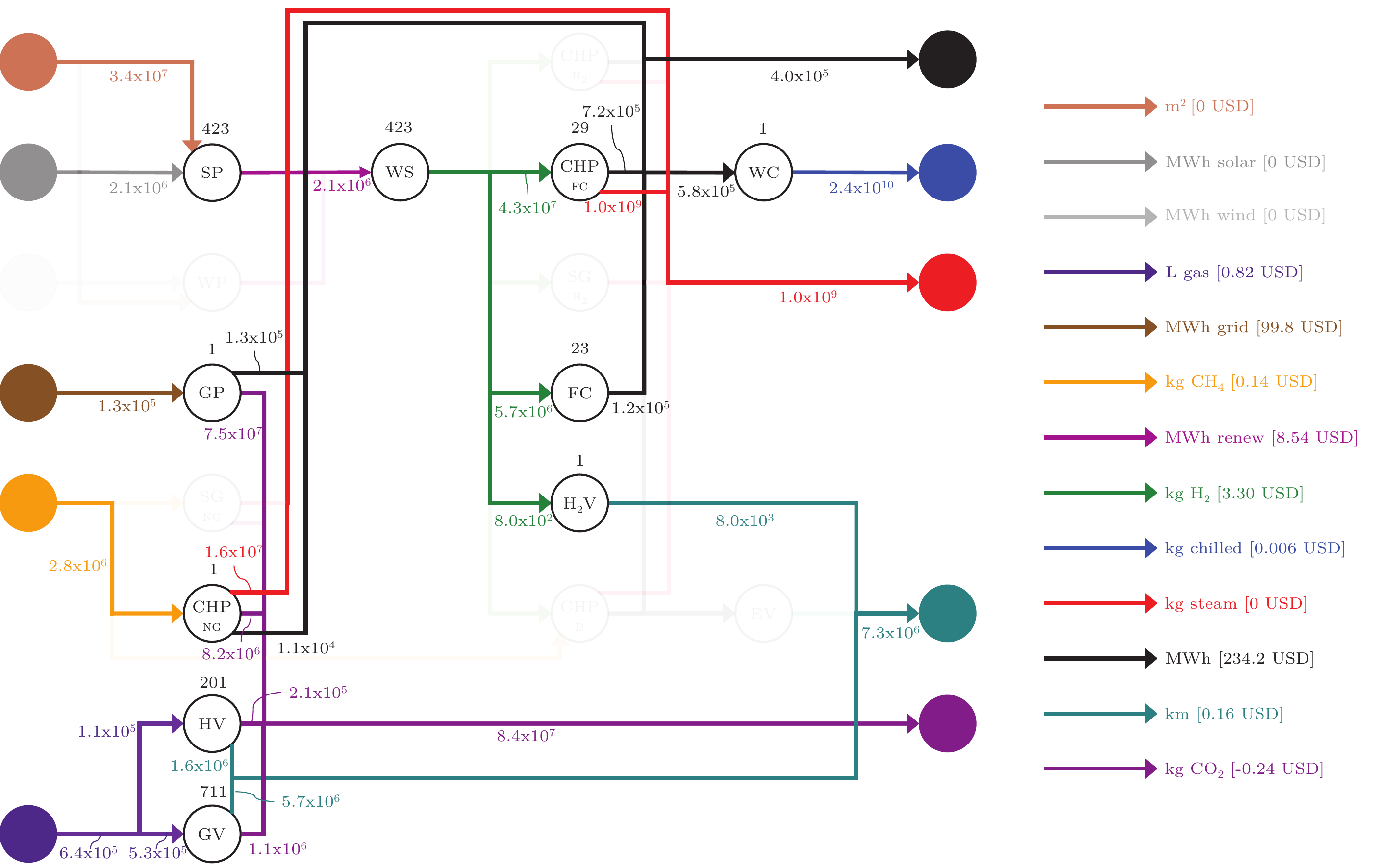}}
\caption{Investment model results at a CO\textsubscript{2} tax of 240 USD per tonne and a budget of 810 million showing the technologies used satisfy the demands of the campus; pathways that are not selected are blurred out. Product prices are displayed in the legend, flows going in and/or out of each stakeholder are reported in the units of the legend, and the number of technology units required are reported above each technology. This pathway makes uses a mix of gasoline, hybrid, and a single hydrogen fuel cell vehicle to meet the distance (km) demand of the university resulting in a high hydrogen price.}
\label{fig:I3}
\end{figure}

The impact of the CO\textsubscript{2} and budget constraint on the hydrogen production and hydrogen price is reported in Figure \ref{fig:IRP}. Hydrogen production is highest in the pathway described in Figure \ref{fig:M2} because it relies on fuel cells to meet the entire electricity use and demand. Hydrogen production is lower when CHP is used, because they are more efficient. The hydrogen price is 0.97 USD per kg; however, it spikes to as high as 3.30 USD per kg when the budget is 810 million USD and the CO\textsubscript{2} tax is 240 USD per tonne. This pathway is displayed in Figure \ref{fig:I3}. The high hydrogen price is caused by the use of a hydrogen vehicle. By participating in the same demand, the hydrogen price becomes tied to the cost of using the most expensive technology used. For this pathway and condition, gasoline vehicles are the most expensive because of the high CO\textsubscript{2} tax and lower efficiency than the hybrid vehicles. This effect on hydrogen price demonstrates the importance of understanding how intermediate products are priced as a result of being used in competing pathways; specifically, the inherent value of intermediate products is inherently tied to the final uses. 
\\

The investment model results demonstrate the impact of the budget constraint on technology selection under different budget and CO\textsubscript{2} tax conditions. These technology selections impact the utility cost for the university and the prices of products. When the budget is low, little can be done to compensate for the CO\textsubscript{2} tax leading to high utility costs and prices. As the budget is increased, technology selections favor the use of solar power production for hydrogen production due to the low operating cost and greatest impact offsetting the CO\textsubscript{2} emissions associated with the use of grid power. Additionally, the investment model results indicate that there are some budgets where potential technologies (e.g., hybrid vehicles) are chosen despite them not being used in the corresponding budget unconstrained case. This suggests that, periodic investments of 100 million USD per year are made for 10 years, may result in a different technology pathway compared to a single investment of 1x10\textsuperscript{9} USD. 
\\

%%%%%%%%%%%%%%%%%%%%%%%%%%%%%%%%%%%%%%%%%%

\section{Conclusions and Future Work}

We have presented an optimization framework for conducting technology pathway analysis. The proposed framework uses a graph abstraction that captures interdependencies between products and technologies and diverse factors that are of interest for high-level analysis such as techno-economic factors, economies of scale, efficiencies, and externalities (e.g., policy). Duality analysis also reveals that the optimization formulation has a natural economic interpretation that reveals how stakeholders participating in the system can be remunerated and that also reveals in inherent value of intermediate products. We have illustrated the use of the framework for studying technology pathways that enable decarbonization of a representative university campus. Here, we show how the framework can model diverse types of products, technologies, and externalities. This analysis identified key technologies and products and demonstrated the importance of understanding product pricing dependencies that can arise as decarbonization occurs. This framework fulfills a key need to be able to quickly screen technologies and understand the complex interactions between technologies in a coherent enabling the exploration of, in particular and exemplified in the case study, emerging and modular technologies. As part of future work, we will explore the use of the proposed framework to decarbonize other systems (such as plastics manufacturing) and to understand the impact of additional technologies (such as ammonia production) and connections with other sectors (such as agriculture). 
\\

\section*{Supporting Information}
Supply, demand, and technology data with references and additional relative figures are available in the supporting information. This information is available free of charge via the Internet at \url{http://pubs.acs.org/}

%%%%%%%%%%%%%%%%%%%%%%%%%%%%%%%%%%%%%%%%%%

\section*{Acknowledgments}

We acknowledge funding from NSF CAREER award CBET-1748516.

%%%%%%%%%%%%%%%%%%%%%%%%%%%%%%%%%%%%%%%%%%

\bibliography{Bib}

%%%%%%%%%%%%%%%%%%%%%%%%%%%%%%%%%%%%%%%%%%

\section*{Abstract Graphic}
\includegraphics[width=1\textwidth]{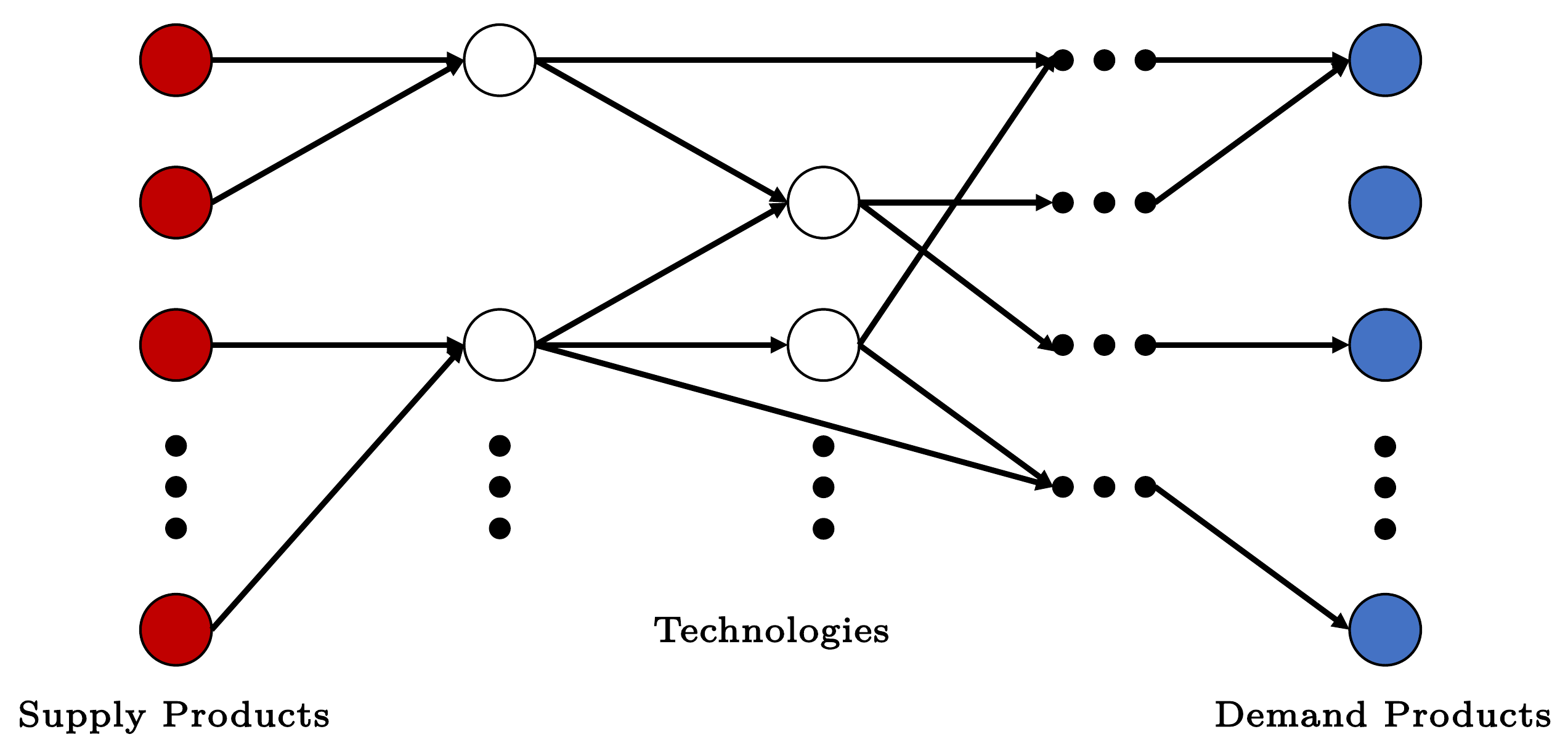}

\end{document}